\documentclass[10pt, reqno]{amsart}
\usepackage{amssymb,mathrsfs,amsmath}
 \usepackage{graphicx}
\usepackage{amscd,amsmath,amsopn,amssymb,amsthm,multicol}
\usepackage{mathtools}
\usepackage{xcolor}
\usepackage{mathrsfs}
\usepackage{color}
\usepackage[matrix, all, 2cell]{xy}
\usepackage{lscape}
\usepackage{slashed}
\usepackage{graphicx}
\usepackage{epsf}
\usepackage{setspace}
\usepackage{upgreek}
\usepackage{textgreek}
\usepackage{tikz}
\usepackage{multirow}
\usepackage{cancel}
\usepackage{soul}
\usepackage{comment}
\usepackage{wasysym}
\usepackage{bbm}
\usepackage{diagbox}
\usepackage{tikz-cd}
\usepackage{bbold}
\usepackage{booktabs}

\advance\oddsidemargin by -1.0cm
\advance\evensidemargin by -1.0cm
\advance\topmargin by -1.0cm
\numberwithin{equation}{section}

\newcommand{\R}{\ensuremath{\mathbb{R}}}
\newcommand{\Z}{\ensuremath{\mathbb{Z}}}
\newcommand{\C}{\ensuremath{\mathbb{C}}}
\newcommand{\Hn}{\ensuremath{\mathbb{H}}}

\renewcommand{\sp}{\mathfrak{sp}}

\newcommand{\la}{\lambda}

\renewcommand{\Re}{\mathbbm{Re}}
\renewcommand{\Im}{\mathbbm{Im}}

\DeclareMathAlphabet{\mathscrbf}{OMS}{mdugm}{b}{n}

\DeclareMathOperator{\SO}{\mathsf{SO}}
\DeclareMathOperator{\Sp}{\mathsf{Sp}}
\DeclareMathOperator{\SU}{\mathsf{SU}}

\DeclareMathOperator{\U}{\mathsf{U}}

\DeclareMathOperator{\Lie}{\mathsf{Lie}}
\DeclareMathAlphabet{\mathpzc}{OT1}{pzc}{m}{it}
\DeclareMathOperator{\Hh}{\mathsf{H}}

\DeclareMathOperator{\E}{\mathsf{E}}
\DeclareMathOperator{\Gl}{\mathsf{GL}}
\DeclareMathOperator{\Sl}{\mathsf{SL}}
\DeclareMathOperator{\Aut}{\mathsf{Aut}}
\DeclareMathOperator{\id}{\mathsf{id}}
\DeclareMathOperator{\Ed}{\mathsf{End}}

\DeclareMathOperator{\Id}{\mathsf{Id}}

\DeclareMathOperator{\ad}{ad}

\DeclareMathOperator{\tr}{tr}
\DeclareMathOperator{\Tr}{\mathsf{Tr}}

\DeclareMathOperator{\Hom}{\mathsf{Hom}}
\DeclareMathOperator{\Ker}{\mathsf{Ker}}

\DeclareMathOperator{\dd}{d}
\newcommand{\fr}{\mathfrak}
\newcommand{\al}{\alpha}

\newcommand{\mc}{\mathcal}

\newcommand{\ep}{\varepsilon}

\newcommand{\om}{\omega}

\usepackage{amsmath}
\DeclareFontFamily{U}{mathx}{}
\DeclareFontShape{U}{mathx}{m}{n}{<-> mathx10}{}
\DeclareSymbolFont{mathx}{U}{mathx}{m}{n}
\DeclareMathAccent{\widehat}{0}{mathx}{"70}
\DeclareMathAccent{\widecheck}{0}{mathx}{"71}

\DeclareMathAlphabet{\mathscrbf}{OMS}{mdugm}{b}{n}

\newcommand{\Gg}{\ensuremath{\mathsf{G}}}
\newcommand{\Lg}{\ensuremath{\mathsf{L}}}
\newcommand{\Kg}{\ensuremath{\mathsf{K}}}

\newcommand{\Rg}{\ensuremath{\mathsf{R}}}
\newcommand{\Ss}{\ensuremath{\mathsf{S}}}

\newcommand{\Trg}{\ensuremath{\mathsf{Tr}}}
\newcommand{\Aff}{\ensuremath{\mathsf{Aff}}}

\newtheorem{theorem}{Theorem}[section]
\newtheorem{lem}[theorem]{Lemma}
\newtheorem{prop}[theorem]{Proposition}
\newtheorem{corol}[theorem]{Corollary}

\theoremstyle{definition}
\newtheorem{defi}[theorem]{Definition}
\newtheorem{example}[theorem]{Example}
 \newtheorem{rem}[theorem]{Remark}
 
\theoremstyle{remark}

\numberwithin{equation}{section}

\def\bd{\begin{defi}}
\def\ed{\end{defi}}
\def\bt{\begin{theorem}}
\def\et{\end{theorem}}
\def\bl{\begin{lem}}
\def\el{\end{lem}}
\def\bp{\begin{prop}}
\def\ep{\end{prop}}
\def\br{\begin{rem}}
\def\er{\end{rem}}
\def\bc{\begin{corol}}
\def\ec{\end{corol}}
\def\bex{\begin{example}}
\def\eex{\end{example}}
\def\pr{\begin{proof}}
\def\pro{\end{proof}}
\def\eqna{\begin{eqnarray*}}
\def\eqnaa{\begin{eqnarray}}
\def\deqna{\end{eqnarray*}}
\def\deqnaa{\end{eqnarray}}

 \definecolor{crew}{rgb}{0.2,0.5,0.2}
\definecolor{red}{rgb}{0.57,0.11,0.15}
\definecolor{cobalt}{rgb}{0.04,0.3,0.85}
\usepackage[backref=page]{hyperref}
\renewcommand*{\backref}[1]{}
\renewcommand*{\backrefalt}[4]{%
	\ifcase #1 %
		\or        (cited on page~#2)%
	\else      (cited on pages~#2) %
	\fi}
\hypersetup{colorlinks=true,linkcolor=cobalt, citecolor=blue,urlcolor=black}

\makeatletter
\def\subsubsection{\@startsection{subsubsection}{3}%
  \z@{.5\linespacing\@plus.7\linespacing}{.3\linespacing}%
  {\normalfont\bfseries}}
\makeatother

\begin{document}

\title[Classification of quaternionic skew-Hermitian symmetric spaces]{Classification of  quaternionic skew-Hermitian symmetric spaces}
\author{Ioannis Chrysikos and Jan Gregorovi\v{c}}
\address{Department of Mathematics and Statistics, 
 	Faculty of Science, Masaryk University,	Kotl\'{a}\v{r}sk\'{a} 2, 611 37 Brno, 
Czech Republic. \\ ORCID ID: 0000-0002-0785-8021 }
\email{chrysikos@math.muni.cz}

\address{Department of Mathematics, Faculty of Science, University of Ostrava, 701 03 Ostrava, Czech Republic, and Institute of Discrete Mathematics and Geometry, TU Vienna, Wiedner Hauptstrasse 8-10/104, 1040 Vienna, Austria. \\ ORCID ID: 0000-0002-0715-7911}\email{ jan.gregorovic@seznam.cz}


\begin{abstract}   
We provide a complete  classification of quaternionic skew-Hermitian symmetric spaces,
namely symmetric spaces that admit a torsion-free $\SO^*(2n)\Sp(1)$-structure for arbitrary $n>1$. 
    Moreover, we prove that any homogeneous   quaternionic skew-Hermitian manifold is necessarily  a symmetric space.
 \end{abstract}
 
 \maketitle 
 



\tableofcontents


\pagestyle{headings}




\section{Introduction}\label{Preliminaries}

 This work is devoted to the classification  of   simply connected symmetric spaces  admitting an invariant $\SO^*(2n)\Sp(1)$-structure,   known as \textsf{quaternionic skew-Hermitian symmetric spaces} (briefly, \textsf{qs-H symmetric spaces}). 
Recall by \cite{CGWPartI}  that an $\SO^*(2n)\Sp(1)$-structure on a $4n$-dimensional manifold $M$ with quaternionic dimension $n>1$ 
is  determined by a pair $(Q, \om)$, where $Q\subset\Ed(TM)$ is an almost quaternionic structure  on $M$ and $\om\in\Lambda^2T^*M$ is a non-degenerate 2-form which is $Q$-Hermitian, usually   referred to as a \textsf{scalar 2-form}. Manifolds admitting such a $\Gg$-structure are called \textsf{almost quaternionic skew-Hermitian manifolds}, abbreviated as \textsf{almost qs-H manifolds}. 
In the torsion-free (i.e., 1-integrable) case   the prefix ``almost'' is omitted  and one speaks of \textsf{quaternionic skew-Hermitian structures}, abbreviated as \textsf{qs-H structures}. 
A  \textsf{quaternionic skew-Hermitian symmetric space}   is then defined as a symmetric space admitting  a quaternionic skew-Hermitian structure $(Q, \om)$ such that each symmetry   is a quaternionic skew-Hermitian automorphism of $(Q, \om)$.   Here, by the term    \textsf{quaternionic skew-Hermitian automorphism} (in short, \textsf{qs-H automorphism}) we mean a diffeomorphism preserving both $Q$ and $\om$. We denote by $\Aut(M, Q, \om)$ the group  of qs-H automorphisms of an (almost) qs-H manifold $(M, Q, \om)$.

In this work we show that these symmetric spaces are the only   homogeneous qs-H manifolds.  

\bt\label{THM1}
Any homogeneous   qs-H manifold  is a qs-H symmetric space.
\et

Moreover, we prove the following classification theorem:
\bt\label{THM2}
Let $M$ be a qs-Hermitian symmetric space. Then, up to covering, $M$ is
one of the following spaces:\\
{\rm(1)} The linear model $[\E\Hh]$ {\rm(}see  \cite[Section 2]{CGWPartI} for details on the linear model{\rm)}. \\
{\rm(2)}  One of the  simple symmetric spaces:
\begin{align*}
M_{1}&=\SO^*(2n+2)/\SO^*(2n)\U(1)\,,\\
M_{2}&=\SU(2+p,q)/(\SU(2)\SU(p,q)\U(1))\,,\\
M_{3}&=\Sl(n+1,\mathbb{H})/(\Gl(1,\mathbb{H})\Sl(n,\mathbb{H}))\,.
\end{align*}
{\rm(3)}  A non-semismple symmetric space $M=\Gg/\Lg$ where $\Gg$ and $\Lg$ have the Levi decompositions 
\[
\Gg=\Sp(p+1, q)\ltimes\Rg\,,\quad \Lg=(\Sp(1)\Sp(p, q))\ltimes\Kg
\]
 for  solvable groups $\Rg$ and $\Kg$, respectively, where the Lie algebras $\fr{r}$ and $\fr{k}=\fr{r}\cap\fr{l}$ of the groups $\Rg$ and $\Kg$ are 
specified according to whether the quaternionic dimension $n$ of $M$ is even or odd:
\begin{itemize}
\item[${\rm(3a)}$] For $n$ even, such symmetric spaces $\Gg/\Lg$ are classified by parameters $p, q$ satisfying $0\leq p, q\leq\frac{n}{2}$  and $0<p+q \leq \frac{n}2$, see Example \ref{clas_exa_1}. 
\item[${\rm(3b)}$] For $n$ odd,  such  symmetric spaces $\Gg/\Lg$ are  classified by parameters $p, q$ satisfying  $0\leq p, q\leq\frac{n-1}{2}$  and $0<p+q \leq \frac{n-1}2$, see Example \ref{clas_exa_2}.
\end{itemize}
For both these classes the quaternionic pseudo-K\"ahler symmetric space $\Sp(p+1, q)/(\Sp(1)\Sp(p, q))$ is the semisimple factor of $\Gg/\Lg$.
\et

Note that the semisimple qs-H symmetric spaces  were previously  classified  in \cite{CGWPartI}, building on the classification results of \cite{AC,G13}. In Section \ref{semisimple_case1} we   re-establish this 
classification by  means of a purely algebraic approach, based  only on the standard theory of symmetric pairs and   extended Dynkin diagrams, see \cite{H78,   CS09} for more details.
 In fact, one can adopt a more general and geometric approach that leads to a classification of all qs-H symmetric spaces, including the semisimple ones. This is based on a construction of local torsion-free $\SO^*(2n)\Sp(1)$-structures, derived by M. Cahen and  L. Schwachh\"ofer    \cite{CahS} (see also \cite[Proposition 6.3]{CGWPartI}). In Section \ref{sec_gen_class} we determine when this construction yields a qs-H symmetric space, see Proposition \ref{general_approach1}.  Then, a combination of Theorem \ref{main_thm_tiLa} and Theorem \ref{iLas_thm}  yields the full classification of
 qs-H symmetric spaces (up to covering), see also Examples \ref{clas_exa_1}-\ref{clas_exa_5} for details.

In Section \ref{section6} we  prove Theorem \ref{THM1}   and, moreover, show   that homogeneous almost qs-H structures do exist. In particular, the construction employed   in this paper yields   homogeneous almost qs-H manifolds  of intrinsic torsion type $\mc{X}_4$  (conformally symplectic structures), see Example \ref{confromal_s} (we refer to  \cite{CGWPartI, CGWPartII} for the intrinsic torsion types of $\SO^*(2n)\Sp(1)$-structures).

  The organization of the paper goes as follows.  
In Section \ref{section2} we review necessary  preliminaries on symmetric spaces. In Section \ref{qsH_tilas} we describe the algebraic structure  equivalent to a  qs-H symmetric space, by introducing the notion of 
\textsf{quaternionic skew-Hermitian transvection involutive Lie algebras}. Using this notion, in Section \ref{semisimple_case1}
we   recover the classification for the semisimple case. Finally, in Sections \ref{sec_gen_class} and \ref{section6} we prove Theorems \ref{THM1} and \ref{THM2}, respectively.
For the sake of brevity,  note that  this  work   does {\it not}  include preliminaries on almost quaternionic skew-Hermitian manifolds. The reader can find all the necessary details in    \cite{CGWPartI, CGWPartII}, see also \cite{CCG}.

\bigskip
\noindent\textbf{Acknowledgements.}
The first author acknowledges  support  from the Czech Grant Agency (project GA24-10887S), and the Horizon 2020 MSCA project CaLIGOLA,   ID 101086123. 
He also thanks the Institute of Mathematics (UZH) at the University of Z\"urich for its hospitality during his academic stay in the Winter 2025 semester.
The second author acknowledges  support by the Austrian Science Fund (FWF) 10.55776/PAT4819724.



\section{The structure of quaternionic skew-Hermitian symmetric spaces}\label{section2}
Our description  begins with a brief summary of symmetric spaces, providing background  relevant  to the discussion that follows.
 There are many equivalent ways to define symmetric spaces.  A review of these approaches is given in  \cite[Section 1.3]{G12} and further  details can be found, for example, in \cite{H78, B00, Kob2}.
 \subsection{Generalities on symmetric spaces}\label{generalities}
\noindent A \textsf{symmetric space} can be viewed as a pair $(M, s)$ where $M$ is a connected smooth manifold  and $s : M\times M\to M$ is a smooth map 
    such that the family 
    \[
    \{s_x : M\to M\}_{x\in M}
    \]
     where $s_{x}(y):=s(x, y)$ for all $x, y\in M$,   is a (smooth) family of involutions of $M$ which satisfies the following properties:
 \begin{itemize}
\item[(S1)] $s_{x}(x)=x$ for all $x\in M$;
\item[(S2)]  $s_{x}\circ s_{y}\circ s_{x}=s_{s_{x}(y)}$ for all $x, y\in M$;
\item[(S3)] For each $x\in M$ there is a neighbourhood $U_{x}\subset M$ where $x$ is the only fixed point of $s_{x}$.
\end{itemize}
Each $s_{x} : M\to M$ is called the \textsf{symmetry} centered at $x$ and by the term  ``involution'' we mean that $s_{x}(s_{x}(y))=y$ for all $x, y\in M$. Note that the third axiom $(S3)$ means that  the point $x$ is an isolated fixed point of $s_x$, for all $x\in M$.  A   \textsf{morphism} between two symmetric spaces $(A, s)$ and  $(B, t)$, with smooth families of involutions $\{s_a\}_{a\in A}$ and $\{t_b\}_{b\in B}$, respectively,  is a smooth map $f : A\to B$ with $f(s_{x}(y))=t_{f(x)}(f(y))$
for any $x, y\in A$.  Such morphisms will be referred to as \textsf{symmetric morphisms}.  The \textsf{automorphism group} of a symmetric space $(M, \{s_x\}_{x\in M})$ is the group of all symmetric automorphisms of $M$, denoted by   $\Aut(M, s)$.  By the second axion $(S2)$ we see that any symmetry $s_x$ centered at $x\in M$  is a symmetric automorphism  of $M$, thus  there is a natural inclusion of the group generated by the symmetries of $(M, s)$ into $\Aut(M, s)$. 
 The smallest subgroup of the group generated by the symmetries   of a symmetric space  $(M, s)$  that acts transitively on  $M$  is  the \textsf{transvection group}, denoted   by $\Trg(M, s):=\langle s_x\circ s_y :  x, y\in M\rangle$. 
 
 For the moment, let us denote this group by  $\Gg=\Trg(M, s)$ and fix a base point $o\in M$. Since $\Gg$ is transitive we have a  homogeneous presentation $M=\Gg/\Lg$, where 
$\Lg=\{g\in \Gg : go=o\}$ is the  stabilizer of   $o\in M$. Note that $\Lg$ is a closed subgroup of $\Gg$.  Moreover, the rule   $\sigma(g):=s_o\circ g\circ s_o$ for all $g\in \Gg$, defines 
an involution $\sigma : \Gg\to \Gg$ of $\Gg$ and the stabilizer $\Lg$ satisfies $(\Gg^{\sigma})^{0}\subset\Lg\subset\Gg^{\sigma}$, 
where $\Gg^{\sigma}$ is the closed subgroup of $\Gg$ consisting of elements fixed by $\sigma$ and $(\Gg^{\sigma})^{0}$ denotes its connected component.  
Conversely, a pair $(\Gg, \sigma)$ of a connected Lie group  $\Gg$ endowed with an (non-trivial) involutive automorphism $\sigma : \Gg\to\Gg$ gives rise to a symmetric space $\Gg/\Lg$, where 
$\Lg$ is any subgroup  with $(\Gg^{\sigma})^{0}\subset\Lg\subset\Gg^{\sigma}$. Often, such triples  $(\Gg, \sigma, \Lg)$ are referred to as  \textsf{symmetric pairs}, see for example \cite[p.~3]{B00}.   
 
  \bd\label{tiLa_Def}
 A finite-dimensional  \textsf{transvection  involutive Lie algebra} (\textsf{tiLa} in short) consists of a finite-dimensional real Lie algebra $\fr{g}$ and an involutive Lie algebra automorphism $\sigma : \fr{g}\to\fr{g}$ such that the  following two condition are satisfied:
\begin{itemize}
\item[(1)] The $\Z_2$-grading of $\fr{g}$, given by
\[
\fr{g}=\fr{l}\oplus\fr{m}\,, \quad \fr{l}:=\{X\in\fr{g} : \sigma(X)=  X\}\,,\quad \fr{m}:=\{X\in\fr{g} : \sigma(X)=- X\},
\]
satisfies the condition $\fr{l}=[\fr{m}, \fr{m}]$.
\item[(2)] No non-zero ideal of $\fr{g}$  is contained in $\fr{l}$.
\end{itemize}
\ed
Given a  tiLa $(\fr{g}, \sigma)$,  the decomposition $\fr{g}=\fr{l}\oplus\fr{m}$ into $\pm 1$-eigenspaces of $\sigma$  is referred to  as the \textsf{canonical decomposition} of $\fr{g}$. Since this decomposition is a $\Z_2$-grading,  the relations  $[\fr{l}, \fr{l}]\subset\fr{l}$ and $[\fr{l}, \fr{m}]\subset\fr{m}$ hold as well. 
Thus $\fr{l}$ is a Lie subalgebra of $\fr{g}$ and  $\fr{m}$ is an  $\fr{l}$-module, which corresponds to a representation 
\[
\rho : \fr{l}\to\Ed(\fr{m})
\]
 defined by  the restriction of the adjoint action.  This is a faithful representation, as a consequence of the second condition in the Definition \ref{tiLa_Def}.  Moreover we mention that  the induced  $\fr{l}$-action on   the endomorphism algebra $\Ed(\fr{m})$ of $\fr{m}$
can be interpreted as  $\rho(X)A:=[\rho(X), A]$, for all  $X\in\fr{l}$ and $A\in\Ed(\fr{m})$. 

 Associated to any symmetric pair  $(\Gg, \sigma, \Lg)$, where $\Gg$ is the transvection group,  there is a unique tiLa  given by    $(\fr{g}:=\Lie(\Gg), \sigma)$.  Here, we maintain the same notation $\sigma$ for the induced involution  at the Lie algebra level, and, as usual, we denote by $\fr{g}=\Lie(\Gg)$ the Lie algebra of $\Gg$.
 \bt\label{FundBij} \textnormal{(\cite{K65, B00})}
 There is a bijective correspondence between tiLas $(\fr{g}, \sigma)$ and symmetric pairs $(\Gg, \sigma, \Lg)$, where  $\Gg$ is the simply connected Lie group corresponding to the Lie algebra $\fr{g}$ and $\Lg$ is the closed connected subgroup of $\Gg$ with  Lie algebra $\fr{l}$. 
Moreover, there is {\color{crew} an!!!} equivalence of categories between isomorphism classes of simply connected symmetric spaces and isomorphism classes of tiLas. 
 \et
Recall that given a smooth connected manifold $M$  with an affine connection $\nabla$,   the group $\Aff(M, \nabla)$ of affine transformations of $(M, \nabla)$  consists of diffeomorphisms $\varphi : M\to M$ preserving $\nabla$, i.e., $\varphi^*\nabla=\nabla$.
An \textsf{affine symmetric space} is a connected manifold $M$ endowed with a  connection $\nabla$ which satisfies the following property for each point $x\in M$: there is an involutive affine transformation $s_{x}\in\Aff(M, \nabla)$ of $(M, \nabla)$ such that $x$ is an isolated point of $s_x$.  Clearly, for each $x\in M$  such an affine transformation $s_x$ is unique,  and it turns out that
an affine symmetric space $(M, \nabla)$ is a symmetric space $(M, s)$, as defined above. Conversely,  any symmetric 
space is an affine symmetric space endowed with its unique \textsf{canonical connection} $\nabla^{0}$   invariant under the symmetries.
Moreover, $\nabla^0$ is torsion-free, complete,  and its curvature tensor $R^{0}$ is $\nabla^{0}$-parallel, see \cite{Kob2}.

\subsection{Quaternionic skew-Hermitian symmetric spaces}
\noindent Symmetric spaces   can  be equipped with various geometric structures
invariant with respect to the symmetries.  For example, a \textsf{pseudo-Riemannian symmetric space} is a symmetric space $(M, \{s_x\}_{x\in M})$  endowed with a pseudo-Riemannian metric $g$ such that each   $s_x$ is an isometry of $g$,
a \textsf{complex symmetric space} is a symmetric space endowed with a complex structure $J$ such that each  $s_x$ is a holomorphic map with respect to $J$, while a \textsf{symplectic symmetric space} is a symmetric space endowed with a symplectic 2-form $\omega$ such that each $s_x$ is a symplectic map with respect to $\omega$.  In this text we are interested in the following type of symmetric spaces.

\bd\label{qsHsym}
A  \textsf{quaternionic skew-Hermitian symmetric space} is a symmetric space admitting  a quaternionic skew-Hermitian structure $(Q, \om)$ such that each $s_x$ is a qs-H automorphism of $(Q, \om)$.   
\ed

Let us recall that, according to  \cite[Theorem 4.12]{CGWPartI}, any almost qs-H manifold $(M^{4n}, Q, \om)$ admits a unique minimal adapted $\SO^*(2n)\Sp(1)$-connection, denoted by $\nabla^{Q, \om}$.    It follows that the  $\SO^*(2n)\Sp(1)$-structure $(Q, \om)$  is 1-integrable if and only if $\nabla^{Q, \om}$ is torsion-free,  see  \cite[Proposition 2.9]{CGWPartII} for more details on the torsion-free case. 
 The reason that we do {\it not} include  the more general almost qs-H structures   in   Definition \ref{qsHsym}   is explained, among other facts,  in the following proposition.

   \bp \label{holonomy_1} \label{canonical_conI} 
Let $(M, s)$ be a symmetric spaced admitting an almost   quaternionic skew-Hermitian symmetric  $(Q, \om)$ invariant under the symmetries. Then,  the  canonical connection $\nabla^0$  has the  following properties:
\begin{itemize}
\item[{\rm (1)}] For all $X, Y, Z\in\fr{X}(M)$ and $x\in M$ we have $\om_x(\nabla^0_{X}Y, Z)=\frac{1}{2}X_{x}(\omega(Y+s_{x_{*}}Y, Z))$, where $s_{x_{*}}=\dd s_{x}$.
\item[{\rm (2)}] The connections $\nabla^{0}$ and $\nabla^{Q, \om}$ coincide, $\nabla^{0}=\nabla^{Q, \om}$ and hence the pair $(Q, \om)$ is necessarily 1-integrable. 
\item[{\rm (3)}]  Any qs-H automorphism $f\in \Aut(M, Q, \om)$ satisfies $f\circ s_{x}=s_{f(x)}\circ f$, for all $x\in M$.  
\item[{\rm (4)}] The linear   holonomy algebra of $\nabla^{0}=\nabla^{Q, \om}$ coincides with the Lie subalgebra $\fr{l}$,  where $\fr{g}=\fr{l}\oplus\fr{m}$ is the canonical decomposition, that is,
\[
\fr{l}\cong{\sf hol}(\nabla^{0})=\{R^{0}(X, Y)=-\ad_{\fr{m}}([X, Y]) : X, Y\in\fr{m}\}\subset\fr{so}^{*}(2n)\oplus\fr{sp}(1)\,.
\]
\end{itemize}
\ep

\pr
We   recall that any qs-H transformation is affine with respect to the $\SO^*(2n)\Sp(1)$-connection $\nabla^{Q, \om}$, see \cite[Proposition 5.8]{CGWPartI}. 
Since $\nabla^{Q, \om}$ is unique and is assumed that is invariant by the symmetries,  this implies that the two connections $\nabla^0$ and $\nabla^{Q, \om}$ coincide, $\nabla^0=\nabla^{Q, \om}$.  The second claim   in (2) follows now since $\nabla^0$ is torsion-free.  Therefore $\om$ is symplectic, and the expression given in (1) can be proved as in   \cite{B98}.  Part (3) follows easily by the fact that   $\nabla^0=\nabla^{Q, \om}$ is invariant under any qs-H  automorphism  $f\in\Aut(M, Q, \om)$. 
Part (4) follows by Theorem 3.2 in \cite[p.~232]{Kob2}.
\pro

By part (3) of Proposition \ref{holonomy_1} it follows that
 \bc
 The category of qs-H symmetric spaces with morphisms the qs-H transformations 
 is a \textsf{full subcategory} of the category of  qs-H manifolds, that is, every qs-H  automorphism of a qs-H symmetric space is a symmetric automorphism. Hence $\Aut(M, s)\cap\Aut(M, Q, \om)=\Aut(M, Q, \om)$.
 \ec

 \section{Quaternionic skew-Hermitian transvection  involutive  Lie algebras}\label{qsH_tilas}
In this  section, we  will establish a bijective correspondence between the isomorphism classes of simply connected quaternionic skew-Hermitian   symmetric spaces   and the isomorphism classes of a following purely
Lie-algebraic object:

\bd\label{qsHiLa}
A \textsf{quaternionic skew-Hermitian transvection involutive Lie algebra} (briefly, \textsf{qs-H tiLa}) is a quadruple $(\fr{g}, \sigma, Q_0, \om_0)$, where:
\begin{itemize}
\item[(1)]    The pair $(\fr{g}, \sigma)$ is a finite-dimensional transvection  involutive Lie algebra.
\item[(2)] The triple   $(\fr{m}, Q_0, \omega_0)$  is an $\fr{l}$-invariant  quaternionic skew-Hermitian vector space, where $\fr{g}=\fr{l}\oplus\fr{m}$  is the associated canonical decomposition of $\fr{g}$. 
\end{itemize}
\ed

\br \label{qsHLa1}
Let $(\fr{g}, \sigma, Q_0, \omega_0)$ be a qs-H tiLa  with $\Z_2$-grading   $\fr{g}=\fr{l}\oplus\fr{m}$, as described above. Next we would  like to remark a few important points:\\ 
(1) Since   $(\fr{g},  \sigma)$ is an involutive Lie algebra  with canonical decomposition $\fr{g}=\fr{l}\oplus\fr{m}$,  we have $B(\fr{l}, \fr{m})=0$, where  $B$ is the Killing form $\fr{g}$. If $B$ is non-degenerate, so is its restrictions  to $\fr{l}$ and $\fr{m}$, respectively,  see \cite[p.~233]{Kob2}.  \\
(2)  By the second condition in  Definition \ref{qsHiLa}, we require  that  both the  geometric objects $Q_0$ and $\om_0$ are  invariant  under the action of the subalgebra $\fr{l}=[\fr{m}, \fr{m}]$.
 In particular, saying that $Q_0\subset\Ed(\fr{m})$ is $\fr{l}$-invariant we mean that $[\rho(X), A]\in Q_0$ for any $X\in\fr{l}$ and $A\in Q_0$.
In analogy,  $\om_0 : \fr{m}\times\fr{m}\to\R$ is a  linear scalar 2-form on $\fr{m}$ with respect to $Q_0$, which is invariant under the action of $\fr{l}$, that is,   $\om_0([X, U], W)+\omega_0(U, [X, W])=0$ for all $X\in\fr{l}$ and $U, V\in\fr{m}$.\\
(3) The extension of $\om_0$ to  the whole $\fr{g}$, which we may still denote by $\om_0$,  is characterized by the relation $X\lrcorner\thinspace\om_0=0$ for all $X\in\fr{l}$.
Moreover, by the $\fr{l}$-invariance of $\om_0$ and the relation $[\fr{m},\fr{m}]=\fr{l}$, it is clear    that  $\om_0$ is a  Chevalley-Eilenberg 2-cocycle  of the trivial representation of $\fr{g}$ on $\R$. This means that
$\om_0$ belongs to the space 
\[
Z^{2}(\fr{g}, \R)=\left\{\al\in C^{2}(\fr{g}, \R)=\Hom(\bigwedge^{2}\fr{g}, \R) : \dd\al=0\right\}
\]
 of Chevalley 2-cocycles  of the trivial representation of $\fr{g}$ on $\R$.
  \er
Let us now define the morphisms in the category of qs-H tiLas. 
\bd\label{isomorphic_iLas}
Two qs-H tiLas $(\fr{g}_{i}, \sigma_i, Q_i, \omega_i)$  $(i=1, 2)$, are said to be \textsf{isomorphic} if there exists a Lie algebra isomorphism $\varphi : \fr{g}_1\to\fr{g}_2$
such that  $\varphi \circ \sigma_1 = \sigma_2 \circ \varphi$,   $\varphi^* Q_2 = Q_1$ and $\varphi^* \omega_2 = \omega_1$. 
\ed

\bt\label{iLas_thm}
There is equivalence of categories between isomorphism classes of qs-H tiLas $(\fr{g}, \sigma, Q_0, \om_0)$ and isomorphism classes  of simply connected qs-H symmetric spaces $(M, Q, \om, s)$.
\et
\pr
  Let $(M=\Gg/\Lg, Q, \om, s)$ be a qs-H symmetric space with $\Gg=\Trg(M, s)$, with  canonical decomposition $\fr{g}=\fr{l}\oplus\fr{m}$, where $\fr{g}=\Lie(\Gg)$ and $\fr{l}=\Lie(\Lg)$, respectively. 
  Let   $\sigma(g)=s_{o}\circ g\circ s_o$  be the involution of $\Gg$ defined  via  the symmetry at the origin $o=e\Lg\in M$.  
  Since $Q\subset\Ed(TM)$ is a $\Gg$-invariant quaternionic structure and $\om$ is a $\Gg$-invariant symplectic 2-form on $M$,  the  corresponding tiLa is given by $(\fr{g}, \sigma, Q_o, \om_o)$, where $\sigma$ is the induced involution on $\fr{g}$ and $(Q_o, \om_o)$ is the evaluation of the pair $(Q, \om)$ at  $o\in M$. We  call  $(\fr{g}, \sigma, Q_o, \om_o)$ the \textsf{qs-H tiLa associated to $(M=\Gg/\Lg, Q, \om, s)$}.  We will now  prove that $(\fr{g}, \sigma, Q_o, \om_o)$ is unique, up to isomorphism in the category of qs-H tiLas.   \\
\noindent {\bf Claim 1:} {\it Let $(M, Q, \om, s)$, $(\hat{M}, \hat{Q}, \hat{\om}, \hat{s})$ be two simply connected  quaternionic skew-Hermitian symmetric spaces  and let  $\psi : M\to \hat{M}$ be a qs-H transformation. 
Let $o\in M$ be a base point in $M$ with associated qs-H tiLa $(\fr{g}, \sigma, Q_{o}, \om_{o})$, and take as base point in $\hat{M}$ the point  $\hat{o}:=\psi(o)$, with associated qs-H tiLa  $(\hat{\fr{g}}, \hat\sigma, \hat Q_{\hat{o}}, \hat\om_{\hat{o}})$. Then the two qs-H tiLas are isomorphic.}\\
\noindent{\it Proof:}
  By definition, and by part (3) in Proposition \ref{holonomy_1},  the qs-H transformation $\psi$ satisfies the relations
   \[
   \psi\circ s_{x}=\hat{s}_{\psi(x)}\circ\psi, \ \ x\in M\,,\quad \psi^*\hat{Q}_{\hat{o}}=Q_{o}\,,\quad \psi^*\hat{\om}_{\hat{o}}=\om_{o}\,.
   \]
  Moreover,     $\psi$ induces a Lie group isomorphism $\Psi : \Gg\to\hat\Gg$ defined by $  \Psi(s_x\circ s_y):=\hat{s}_{\psi(x)}\circ\hat{s}_{\psi(y)}$, 
   for all $ x, y\in M$. Hence $\Psi$ maps the generators of $\Gg=\Tr(M, s)$ to the generators of $\hat\Gg=\Tr(\hat{M}, \hat{s})$.  This satisfies $\Psi\circ \sigma=\hat{\sigma}\circ\Psi$,   and it follows that the differential $\dd\Psi_{e} : \fr{g}\to \hat{\fr{g}}$  is the desired  isomorphism between the associated qs-H tiLas $(\fr{g}, \sigma, Q_o, \om_o)$ and $(\hat{\fr{g}}, \hat\sigma, \hat Q_{\hat{o}}, \hat\om_{\hat{o}})$. \hfill $\square$
   
\noindent   This proves the one direction of the statement in the theorem. 
  Let us now prove the converse.  Let $(\fr{g}, \sigma, Q_{0}, \om_{0})$  be a qs-H tiLa and let  $\fr{g}=\fr{l}\oplus\fr{m}$ be the corresponding $\Z_2$-grading of $\fr{g}$.  Let $\Gg$ be the  connected and simply connected Lie group corresponding 
to the Lie algebra $\fr{g}$, and let $\sigma : \Gg\to\Gg$ be the  analytic automorphism of $\Gg$ generated by $\sigma$. This is an involutive automorphism of $\Gg$ and since $\Gg$ is simply connected,  the fixed point set $\Gg^{\sigma}$ is connected (see \cite[p.~293]{K65}). Hence its is clear that  $\Lg=\Gg^{\sigma}$ is the analytic subgroup generated by the Lie subalgebra $\fr{l}$ and moreover that  the coset $M:=\Gg/\Lg$ 
is a simply connected effective  symmetric space, see also Theorem \ref{FundBij}. 

\noindent Let $o=e\Lg$ be the origin in $M=\Gg/\Lg$. 
The pair  $(Q_{0}, \om_{0})$ is by definition an $\fr{l}$-invariant linear qs-H structure on $\fr{m}$, which is also $\Lg$-invariant  (since $\Lg$ is connected).  Moreover, the pair  $(Q_{0}, \om_{0})$  is invariant under the symmetry $s_{o}$.
Thus it descends to a $\Gg$-invariant qs-H structure $(Q, \om)$ on $\Gg/\Lg$, with $Q_{o}=Q_{0}$ and $\om_{o}=\om_{0}$, respectively, which is invariant under the symmetries. It remains to prove that  $(M=\Gg/\Lg,  Q, \om, s)$ is unique  up to isomorphism. \\
{\bf Claim 2:} {\it Let  $(\fr{g}, \sigma, Q_{0}, \om_{0})$  and $(\hat{\fr{g}}, \hat\sigma, \hat Q_{0}, \hat\om_{0})$ be two isomorphic qs-H tiLas. Then the corresponding  simply connected qs-H symmetric spaces
$(M=\Gg/\Lg,  Q, \om, s)$ and $(\hat{M}=\hat{\Gg}/\hat{\Lg}, \hat{Q}, \hat{\om}, \hat{s})$ 
are isomorphic.}\\
\noindent{\it Proof:} Let 
us denote the  isomorphism between the qs-H tiLas by $\varphi : \fr{g}\to\hat{\fr{g}}$. By definition,  $\varphi\circ\sigma=\hat\sigma\circ\varphi$ and $\varphi^*\hat{Q}_{0}=Q_{0}$, $\varphi^*\hat\om_{0}=\om_{0}$. Moreover, $\varphi$ preserves the $\Z_2$-grading and we have $\varphi(\fr{l})=\hat{\fr{l}}$ and $\varphi(\fr{m})=\hat{\fr{m}}$.
Let $\Gg$ and $\hat\Gg$ be the connected, simply connected Lie groups with Lie
algebras $\fr{g}$ and $\hat{\fr g}$, respectively, and let $\sigma,\hat\sigma$ be the corresponding involutions at the Lie group level.  The isomorphism $\varphi$ integrates to a unique 
Lie group isomorphism  $\Phi : \Gg\to\hat{\Gg}$, 
 satisfying $\Phi\circ\sigma=\hat\sigma\circ\Phi$.
Therefore $\Phi(\Gg^\sigma)=\hat\Gg^{\hat\sigma}$, and in particular $\Phi(\Lg)=\hat\Lg$.
Hence we obtain a diffeomorphism $f : \Gg/\Lg\to\hat{\Gg}/\hat{\Lg}$ with $f(g\Lg)=\Phi(g)\hat{\Lg}$, $g\in\Gg$. Its differential at the origin $o\in\Gg/\Lg$ satisfies $\dd f_{o}=\varphi|_{\fr{m}}$
and since $\varphi$ is an isomorphism of qs-H tiLas, it follows that $f^*\hat Q=Q$ and $f^*\hat\om=\om$, where $(Q, \om)$ is the $\Gg$-invariant qs-H structure corresponding to 
$(Q_{0}, \om_{0})$ and $(\hat{Q}, \hat{\om})$ is the $\hat{\Gg}$-invariant qs-H structure corresponding to $(\hat Q_{0}, \hat\om_{0})$.
To complete the proof one can   use  part (4) of Proposition \ref{holonomy_1}, which implies that $f : \Gg/\Lg\to\hat{\Gg}/\hat{\Lg}$  satisfies  $f\circ s_{g\Lg} = \hat s_{f(g\Lg)}\circ f$
(this  also follows from the relation $\Phi\circ\sigma=\hat\sigma\circ\Phi$).
\pro

 \section{Classification in the semisimple case}\label{semisimple_case1}
 Theorem \ref{iLas_thm} implies that the classification of simply connected qs-H symmetric spaces   is equivalent to the classification of qs-H tiLas.  Next we will rely on this fact to classify all  simply connected qs-H symmetric spaces  (up to covering).
 To this end we will first treat the semisimple case,  hence  it is sufficient to focus on quaternionic skew-Hermitian tiLas $(\fr{g}, \sigma, Q_0, \om_0)$ with $\fr{g}$ semisimple.

\subsection{Structure of semisimple quaternionic skew-Hermitian tiLas}
 \noindent      Let   $(\fr{g}, \sigma, Q_0, \om_0)$ be a quaternionic skew-Hermitian tiLa  with canonical decomposition $\fr{g}=\fr{l}\oplus\fr{m}$.  
   If $\fr{g}$ is semisimple, then  its Killing form $B$ is non-degenerate,  and as we mentioned above, the  same holds for the restrictions   of $B$ to $\fr{l}$  and $\fr{m}$, denoted by $B_{\fr{l}}:=B|_{\fr{l}\times\fr{l}}$ and   $B_{\fr{m}}=B|_{\fr{m}\times\fr{m}}$, respectively.  The   proposition below  is a structural result regarding the  qs-H tiLas when   $\fr{g}$ is simple.

\bp\label{iLa's_center} Let $(\fr{g}, \sigma, Q_0, \om_0)$  be a qs-H tiLa with $\fr{g}$ simple and  canonical decomposition $\fr{g}=\fr{l}\oplus\fr{m}$.
  Then, there is a non-zero element  $Z_0\in\fr{g}$  belonging to the center $Z(\fr{l})$ of $\fr{l}$.
  \ep
  \pr
 Under our assumptions, both $B_{\fr{m}}$ and $\om_0$ are non-degenerate, and one can  naturally define an endomorphism $I\in\Ed(\fr{m})$ by the relation 
  \[
  B_{\fr{m}}(X, Y):=\om_0(X, IY)\,,\quad X, Y\in\fr{m}\,.
  \]
  \noindent{\bf Claim 1:}  {\it The endomorphism $I$ commutes with the Lie subalgebra $\fr{l}$.} \\
  {\it Proof.} Recall that $\fr{l}$ acts on $\fr{m}$ via the faithful representation $\rho : \fr{l}\to\Ed(\fr{m})$ defined by the adjoint action, $\rho(A)X=[A, X]$ for all $A\in\fr{l}$ and $X\in\fr{m}$. Saying that $I$ commutes with $\fr{l}$ we mean that $I\circ\rho(A)=\rho(A)\circ I$, that is
  \begin{equation}\label{A-invariance}
  I[A, X]=[A, IX]\,,\quad \forall A\in\fr{l}, \ X\in\fr{m}\,.
  \end{equation}
 To prove this, note that the restriction $B_{\fr{m}}$ is   $\ad_{\fr{l}}$-invariant, i.e.,   $B_{\fr{m}}([A, X], Y)=-B_{\fr{m}}(X, [A, Y])$, for all $A\in\fr{l}$, and  $X, Y\in\fr{m}$.  By the definition of $I$ we can equivalently express this as $\om_0([A, X], IY)=-\om_0(X, I[A, Y]) $, that is, 
  \[
  \om_0([A, X], IY)+\om_0(X, I[A, Y])=0\,,\quad\quad  A\in\fr{l}, \ X, Y\in\fr{m}\,.\quad\quad(\dag)
  \]
  Moreover, by the $\ad_{\fr{l}}$-invariance of $\om_0$ we have $\om_0([A, X], Z)+\om_0(X, [A, Z])=0$ for all $A\in\fr{l}$, and $X, Z\in\fr{m}$. Thus, for $Z=IY$ we get 
  \[
  \om_0([A, X], IY)+\om_0(X, [A, IY])=0\,, \quad\quad  A\in\fr{l}, \ X, Y\in\fr{m}\,.\quad\quad(\ddag)
  \]
  By $(\dag)$ and $(\ddag)$ we obtain $\om_0(X, I[A, Y])=\om_0(X, [A, IY])$ which  implies (\ref{A-invariance})  by the non-degeneracy of $\om_0$. This proves Claim 1.\hfill$\square$\\
  \noindent  Let us now extend the endomorphism $I$ to the whole $\fr{g}$ by setting $I(\fr{l})=0$.  By  (\ref{A-invariance}) and the inclusions $[\fr{l}, \fr{m}]\subset\fr{m}$, $[\fr{m}, \fr{m}]\subset\fr{l}$, we see that
  \[
  I[X+A,Y+C]=I[X,C]+I[A,Y]=[IX,C]+[A,IY]\,,
  \]
  and moreover   $[I(X+A),Y+C]+[X+A,I(Y+C)]=[IX,C]+[A,IY]+[IX,Y]+[X,IY]$, 
  for all $A, C\in\fr{l}$ and  $X, Y\in\fr{m}$. Therefore, if the relation $[IX,Y]+[X,IY]=0$ holds,  the previous two relations will imply that 
  \[
   I[X+A,Y+C]=[I(X+A),Y+C]+[X+A,I(Y+C)]\,,
  \]
which means that  $I$ should be a derivation of $\fr{g}$.  Let us prove that indeed this is the case.
First we will show that:\\
\noindent {\bf Claim 2:}  {\it The endomorphism $I : \fr{m}\to\fr{m}$ is skew-symmetric with respect to $B_{\fr{m}}$.} \\
\noindent {\it Proof.} Indeed, since $B_{\fr{m}}$ is symmetric and $\om_0$ skew-symmetric, we immediately get   
\begin{equation}\label{omIX}
\om_0(X, IY)=-\om_0(IX, Y)\,,\quad X, Y\in\fr{m}\,.
\end{equation}
Then   $B_{\fr{m}}(IX, Y)\overset{{\rm def}}{=}\om_0(IX, IY)$, $B_{\fr{m}}(X, IY)\overset{{\rm def}}{=}\om_0(X, I(IY))\overset{(\ref{omIX})}{=}-\om_0(IX, IY)$, 
and by comparing these two relations we deduce that $B_{\fr{m}}(IX, Y)+B_{\fr{m}}(X, IY)=0$ for all $X, Y\in\fr{m}$. \hfill$\square$\\
\noindent Let now $\Gg/\Lg$ be the (simply connected) effective coset associated to $(\fr{g}, \sigma, Q_0, \om_0)$ as constructed in the proof Theorem \ref{iLas_thm},  and let $\nabla^0$ be the    canonical connection associated to the canonical decomposition $\fr{g}=\fr{l}\oplus\fr{m}$, with curvature tensor $R^{0}$. Recall that 
$R^{0}(X, Y)Z=-[[X, Y], Z]$  for all $X, Y, Z\in \fr{m}\cong T_{o}\Gg/\Lg$.  
By the properties of the curvature tensor $R^{0}$  and by Claims 1 and 2,  we see that
\begin{align*}
B_{\fr{m}}([[IX, Y], U], V)&=-B_{\fr{m}}(R^{0}(IX, Y)U, V)=-B_{\fr{m}}(R^{0}(U, V)IX, Y)=B_{\fr{m}}([[U, V], IX], Y)\\
&\overset{{\rm Claim \ 1}}{=} B_{\fr{m}}(I[[U, V], X], Y)\overset{{\rm Claim \ 2}}{=}-B_{\fr{m}}([[U, V], X], IY)\\
&=B_{\fr{m}}(R^{0}(U, V)X, IY)= B_{\fr{m}}(R^{0}(X, IY)U, V)=-B_{\fr{m}}([[X, IY],U], V)
\end{align*}
for all  $X, Y, U, V\in\fr{m}$, which implies that  $[IX, Y]+[X, IY]=0$ for all $X, Y\in\fr{m}$.   Therefore,  $I\in\Ed(\fr{m})$  is a derivation of $\fr{g}$ and since $\fr{g}$  is  simple,  there exists an
element $Z_0\in\fr{g}$ such that $I=\ad(Z_0)$.  By Claim 1 we then deduce that $Z_{0}\in Z(\fr{l})$.  This completes the proof.
\pro

Suppose now that  $(\fr{g},\sigma,Q_0, \om_0)$ is \textsf{semisimple qs-H tiLa}, that is, a qs-H tiLa with semisimple Lie algebra $\fr{g}$.   Let   
\[
\fr{g}=\fr{g}_1\oplus\cdots\oplus\fr{g}_s
\]
be the decomposition of $\fr{g}$ into simple ideals, that is, either each $\fr{g}_{i}$ is a simple Lie algebra and  $\sigma(\fr{g}_i)=\fr{g}_i$, or $\fr{g}_{i}=\fr{g}'_{i}\oplus\fr{g}'_{i}$, where $\fr{g}'_i$ is simple and $\sigma$ satisfies  $\sigma((X, Y))=(-Y, X)$, for all $X, Y\in\fr{g}'_{i}$. 
  Then we can decompose $\sigma$ as  $\sigma = \sigma_1 \oplus \cdots \oplus \sigma_s$,
  with  $\sigma_i := \sigma|_{\fr{g}_i} : \fr{g}_i \to \fr{g}_i$ for each $i$, satisfying $\sigma_i^2 = \Id_{\fr{g}_i}$, for all $i\in\{1, \ldots, s\}$. We can also consider the $\Z_2$-decompositions $\fr{g}_{i}=\fr{l}_{i}\oplus\fr{m}_{i}$ with $\fr{l}_{i}:=\fr{g}_{i}^{\sigma_i}$ and $\fr{m}_{i}:=\fr{g}_{i}^{-\sigma_i}$.   Obviously, we have $\fr{l}_{i}=\fr{l}\cap\fr{g}_i$ and $\fr{m}_i=\fr{m}\cap\fr{g}_{i}$ for all $i$, where $\fr{l}=\fr{g}^{\sigma}$ and $\fr{m}=\fr{g}^{-\sigma}$, respectively. Then we can write
\begin{equation}\label{semi_iLa}
\fr{g}=\bigoplus_{i=1}^{s}\fr{g}_{i}
=\bigoplus_{i=1}^{s}(\fr{l}_i\oplus\fr{m}_i)=\bigoplus_{i=1}^{s}\fr{l}_i\oplus\bigoplus_{i=1}^{s}\fr{m}_i=\fr{l}\oplus\fr{m}\,,
\end{equation}
with $\fr{l}=\bigoplus_{i=1}^{s}\fr{l}_i$ and $\fr{m}=\bigoplus_{i=1}^{s}\fr{m}_i$, respectively.

 \bp\label{La_simple1}\label{ILa_simple1}
For a semisimple  qs-H tiLa $(\fr{g},\sigma,Q_0,\omega_0)$ the Lie algebra $\fr{g}$ is necessarily simple.
 \ep
 \pr
Since $Q_0$ is a linear quaternionic structure we may choose an admissible basis of $Q_0$ and identify the simple Lie algebra $\fr{sp}(1)$ with $Q_0$, $\fr{sp}(1)\cong Q_0$. The specific choice of this admissible basis will not affect the argument.  
 Assume, to the contrary, that  $s>1$.  Then we have that\\
{\bf Claim 1:} {\it $Q_0$ cannot commute with $\fr{l}_{i}$ for all $i$,  and in particular $Q_0$   does not commute with $\fr{l}$.}\\
\noindent{\it Proof.}
Suppose, to the contrary,  that  $Q_0\subset\Ed(\fr{m})$  commutes with  each Lie subalgebra $\fr{l}_i$ for  $i\in\{1, \ldots, s\}$. Then $Q_0$ preserves $\fr{m}_{i}$, which means that $A(\fr{m}_{i})\subset\fr{m}_{i}$ for all $A\in Q_0$ and all $i\in\{1, \ldots, s\}$.
Hence in this case  each $\fr{l}_{i}$-module $\fr{m}_{i}$ should be of quaternionic type (a quaternionic representation),   which contradicts \cite[Corollary V.1.13]{B00} (since each $\fr{g}_{i}$ is a simple ideal).   Hence $Q_0$ does not commute with  all $\fr{l}_{i}$, and the claim for $\fr{l}$ follows.  \hfill $\square$\\
\noindent
Now,  by the definition of a quaternionic skew-Hermitian tiLa, $Q_0$ is    $\fr{l}$-invariant which means that  $[X, A]\in Q_0$ for any $X\in\fr{l}$ and $A\in Q_0$ (see Remark \ref{qsHLa1}).    Since $Q_0\cong\fr{sp}(1)$ and $\fr{sp}(1)$ is simple, this means that  the adjoint action defines a Lie algebra homomorphism
\[
\lambda : \fr{l}\to {\sf Der}(Q_0)=\ad(Q_0)\cong\fr{sp}(1)\,,\quad \lambda(X)A=[X, A]\,.
\]
  Restricting $\lambda$ to each $\fr{l}_i$ we obtain representations $\la_i : \fr{l}_i\to \fr{sp}(1)$. Since $Q_0$ is a simple $\fr{l}_{i}$-module,  the image of each $\fr{l}_i$ should be either trivial or  isomorphic to $\fr{u}(1)$ or to  $\fr{sp}(1)$ itself,   i.e.,  $\la_i(\fr{l}_i)=\{0\}$ or $\la_i(\fr{l}_i)=\fr{u}(1)$ or $\la_i(\fr{l}_i)=\fr{sp}(1)$.  Because $\fr{l}$ acts non-trivially on $Q_0$ and $\fr{sp}(1)$ has no non-trivial commuting subalgebras,  at least one $\fr{l}_i$ should act non-trivially on $Q_0$,  hence we have either $\la_i(\fr{l}_i)=\fr{u}(1)$ for some $i$, or  $\la_i(\fr{l}_i)=\fr{sp}(1)$ for some $i$.  Then since  $[\fr{l}_i, \fr{l}_j]=0$ for all $i\neq j$  it is easy to see that  $\la_j(\fr{l}_j)=0$ for all $j\neq i$, which implies that    $[\fr{l}_j, Q_0]=0$. This gives a contradiction according to Claim 1.  Hence $s=1$ and $\fr{g}=\fr{g}_i$. The case  $\fr{g}_{i}=\fr{g}'_{i}\oplus\fr{g}'_{i}$ with $\fr{g}'_{i}$ simple is not possible,  because the adjoint representation of a simple Lie algebra is never of quaternionic type.
  This now implies that $\fr{g}$ is simple.
   \pro

\subsection{Classification of semisimple qs-H symmetric spaces}
\noindent  We now turn our attention to the classification of semisimple qs-H symmetric spaces.   

 \bt\label{ClasThm} \textnormal{(\cite{CGWPartI})} A   semisimple quaternionic skew-Hermitian symmetric space $(M=\Gg/\Lg, Q, \om, s)$ is, up to covering, one of the following symmetric spaces: 
\begin{align*}
M_{1}&=\SO^*(2n+2)/\SO^*(2n)\U(1)\,,\\
M_{2}&=\SU(2+p,q)/(\SU(2)\SU(p,q)\U(1))\,,\\
M_{3}&=\Sl(n+1,\mathbb{H})/(\Gl(1,\mathbb{H})\Sl(n,\mathbb{H}))\,.
\end{align*}
Moreover, there are the following additional geometric structures provided by the central part of $\Lg$:
\begin{itemize}
\item $M_1$ admits a distinguished $Q$-compatible invariant  complex structure $J$ and the corresponding invariant metric $g_{J}=\om\circ J$ is pseudo-K\"ahler.
\item $M_2$ admits a distinguished invariant complex structure $I$ commuting with $Q$ and the corresponding invariant metric $g_{I}=\om\circ I$ is quaternionic pseudo-K\"ahler. For $q=0$, they space $M_2$ is compact and $g_{I}$ is  quaternionic K\"ahler.  
\item $M_3$ admits a distinguished invariant para-complex structure $I$ commuting with $Q$ and the corresponding invariant metric $g_{I}=\om\circ I$ is quaternionic pseudo-K\"ahler of split signature. 

\end{itemize}
\et
\pr
As we mentioned above, to obtain the classification of semisimple qs-H symmetric spaces it  suffices to classify the semisimple  
quaternionic skew-Hermitian tiLas and our approach for this relies on the standard classification of  symmetric pairs $(\fr{g}, \sigma)$ for $\fr{g}$ semisimple. This can be carried out using extended Dynkin diagrams, following the approach of the classical monograph \cite{H78}. 

 Let $(\fr{g}, \sigma, Q_0, \om_0)$ be a  semisimple   qs-H tiLa, with canonical decomposition $\fr{g}=\fr{l}\oplus\fr{m}$. 
By  Proposition \ref{La_simple1}  the Lie algebra $\fr{g}$ is simple and by Proposition \ref{iLa's_center}   the center of the Lie subalgebra $\fr{l}$ is necessarily non-trivial  $Z(\fr{l})\neq\{0\}$. 
 This fact simplifies our examination, as it reduces the possible extended Dynkin diagrams  that one need to consider, only  to those   in Table 1  of \cite[p.~503]{H78}.  
 
The procedure depends on the intersection of $Q_0\cong\fr{sp}(1)$   with $\rho(\fr{l})$, which we denote by
 \[
 \fr{s}:=\rho(\fr{l})\cap\fr{sp}(1)\,, 
 \]
where $\rho : \fr{l}\to\Ed(\fr{m})$ is the faithful representation associated to the quadruple $(\fr{g}, \sigma, Q_0, \om_0)$.    We know that this intersection cannot be  trivial; if $\fr{s}=\{0\}$ then  $\rho$ should be a quaternionic representation, which is impossible   by \cite[Corollary V.1.13]{B00}.  Actually, since $\fr{sp}(1)$ does not contain any 2-dimensional proper Lie subalgebra   the  intersection $\fr{s}$ can be either 1-dimensional or 3-dimensional. 

\noindent {\bf Case A ($\dim\fr{s}=1$).}  If  $\dim\fr{s}=1$, we are looking in Table 1  of \cite[p.~503]{H78} for a case where, when we exclude the center of $\fr{l}$, then the representation is quaternionic. This is a classical problem that can be easily resolved, because it is easy to find out the representation in question from Table 1 in \cite[p.~503]{H78} and  next use the tables in \cite[pp.~612-615]{CS09}    to determine, whether the representation is quaternionic or not.  The recipe is  as follows:\\
\noindent $\bullet$ Consider the two nodes in the extended Dynkin diagram determining $\sigma$ and  assign the fundamental weight $1$, or $2$, according to the recipe in \cite[Section 3.2.18]{CS09}  to the attached nodes to describe the representation of the  Lie subalgebra   $\fr{l}$    on $\fr{m}$.   Note that $\fr{l}$ is defined by the remaining nodes in the extended Dynkin diagram.\\
\noindent $\bullet$ Next, use the index information from the tables in \cite{CS09} to determine  whether  the total index is $-1$,  and so    the representation is quaternionic.\\
\noindent Such a  procedure finally provides only one case  that is admissible,   namely this with extended Dynkin diagram
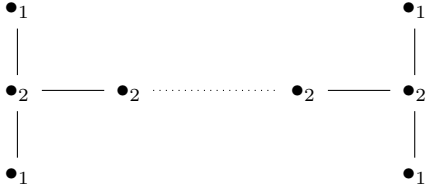
\[\begin{tikzcd}
 {\bullet_1} &&&& {\bullet_1} \\
 {\bullet_2} & {\bullet_2} && {\bullet_2} & {\bullet_2} \\
 {\bullet_1} &&&& {\bullet_1}
 \arrow[no head, from=1-1, to=2-1]
 \arrow[no head, from=1-5, to=2-5]
 \arrow[no head, from=2-1, to=3-1]
 \arrow[no head, from=2-2, to=2-1]
 \arrow[dotted, no head, from=2-2, to=2-4]
 \arrow[no head, from=2-4, to=2-5]
 \arrow[no head, from=2-5, to=3-5]
\end{tikzcd}\]
Hence there is only admissible  symmetric pair with $\dim\fr{s}=1$,  namely  
\[
 (\fr{so}^*(2n+2), \fr{so}^*(2n)\oplus\fr{u}(1))\,,
 \]
which corresponds to  the symmetric space  $M_1$.   \\ 
{\bf Case B ($\dim\fr{s}=3$).}  If $\dim\fr{s}=3$, then we are looking in Table 1  of \cite[p.~503]{H78}   for a case where, after excluding an $\fr{sp}(1)$-factor from $\fr{l}$,   the resulting representation is quaternionic. 
 In a manner similar to the previous case, we compare  this table with  the tables in \cite{CS09} to determine, whether the representation is quaternionic.  In this case the comparison provides that only the extended Dynkin diagrams of $A$-type are admissible, having the form
  
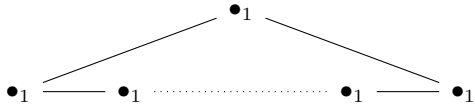
\[\begin{tikzcd}
 && {\bullet_1} \\
 {\bullet_1} & {\bullet_1} && {\bullet_1} & {\bullet_1}
 \arrow[no head, from=1-3, to=2-5]
 \arrow[no head, from=2-1, to=1-3]
 \arrow[no head, from=2-2, to=2-1]
 \arrow[dotted, no head, from=2-2, to=2-4]
 \arrow[no head, from=2-4, to=2-5]
\end{tikzcd}\]
We deduce that there are two admissible symmetric pairs with $\dim\fr{s}=3$, namely
  \[
  (\fr{su}(2+p, q), \fr{su}(2)\oplus\fr{su}(p, q)\oplus\fr{u}(1))\,,\quad (\fr{sl}(n+1, \Hn), \fr{gl}(1, \Hn)\oplus\fr{sl}(n, \Hn))
  \] 
  These pairs correspond to the symmetric spaces $M_2$ and $M_3$ respectively. For the characterization of the additional geometric structure on each of the three symmetric spaces we refer to \cite{CCG}, see also \cite{Pont} for the case of $M_1$. This completes our proof. 
 \pro

\br\label{remark_c}
 Let $(M=\Gg/\Lg, Q, \om, s)$ be a semisimple qs-H symmetric space. Then the center of $\Lg$ incudes either an invariant complex or paracomplex structure on $M$. Therefore, it follows that $\fr{l}$ is actually the centralizer of $Z_{0}\in Z(\fr{l})$ in $\fr{g}$, i.e., $\fr{l}\cong C_{\fr{g}}(Z_0)$.
\er

\section{Classification in the non-semisimple case}\label{sec_gen_class}
The classification of general qs-H symmetric spaces relies on a geometric 
 construction described in \cite{CahS},  which provides local qs-H manifolds (see also  \cite[Section 6.2]{CGWPartI}).  
 As we will see below, this method  yields the semisimple qs-H symmetric spaces, as well. 

The starting point is the \textsf{contact grading} of the Lie algebra $\fr{so}^*(2n+4)$,
given as follows:
 We   represent   $\fr{so}^*(2n+4)$ by  matrices of the form
\begin{align*}
\begin{pmatrix}
a & Y^*\mathbb{j} & d\\
X & A & Y\\
 c &-X^*\mathbb{j} & -a^*
\end{pmatrix},
\end{align*}
where   $a\in\Hn$,  $X, Y\in \Hn^n$, $c, d\in \R$ and $A\in \fr{so}^*(2n)$. 
Here,  $\mathbb{j}$ denotes the matrix of the   scalar 2-form on $\Hn^n$   and $\al^*$ denotes the conjugate transpose   of $\al$.   Then the contact grading is described by
\[
\fr{so}^{*}(2n+4)=\frak{so}^*(2n+4)_{-2}\oplus\fr{so}^{*}(2n+4)_{-1}\oplus \fr{so}^*(2n+4)_{0}\oplus \fr{so}^*(2n+4)_{1}\oplus\fr{so}^*(2n+4)_{2}
\]
with
\begin{eqnarray*}
\fr{so}^{*}(2n+4)_{-2}&:=&\left\{
\begin{pmatrix}
0 & 0 & 0\\
0 & 0 & 0\\
c & 0 & 0
\end{pmatrix} : c\in\R\right\}\cong\R,\\
\fr{so}^{*}(2n+4)_{-1}&:=&
\left\{
\begin{pmatrix}
0 & 0 & 0\\
X & 0 & 0\\
0 & -X^*\mathbb{j} & 0
\end{pmatrix} : X\in\Hn^n\right\}\cong\Hn^n,\\
\fr{so}^{*}(2n+4)_{0}&:=&
\left\{
\begin{pmatrix}
a & 0 & 0\\
0 & A & 0\\
0 & 0 & -a^*
\end{pmatrix} : a\in\Hn,\; A\in\fr{so}^*(2n)\right\}\cong \Hn\oplus\fr{so}^*(2n),\\
\fr{so}^{*}(2n+4)_{1}&:=&
\left\{
\begin{pmatrix}
0 & Y^*\mathbb{j} & 0\\
0 & 0 & Y\\
0 & 0 & 0
\end{pmatrix} : Y\in\Hn^n \right\}\cong\Hn^n,\\
\fr{so}^{*}(2n+4)_{2}&:=&
\left\{
\begin{pmatrix}
0 & 0 & d\\
0 & 0 & 0\\
0 & 0 & 0
\end{pmatrix} : d\in\R\right\}\cong\R\,.
\end{eqnarray*}
Note that $\fr{g}_{\pm2}$ are 1-dimensional and the bilinear form (Levi-bracket)  $[  \cdot , \cdot  ] : \fr{g}_{-1}\times\fr{g}_{-1}\to\fr{g}_{-2} $ is a $\fr{g}_{-2}$-valued scalar 2-form, and hence non-degenerate. 
Set  $e:=\begin{pmatrix}
0&0&1\\
0&0&0\\
0&0&0
\end{pmatrix}\in\fr{g}_{2}$ and $f:=
\begin{pmatrix}
0&0&0\\
0&0&0\\
1&0&0
\end{pmatrix}\in\fr{g}_{-2}$.  The grading element $h\in\fr{g}_{0}$ is given by 
\[
h=[e, f]=\begin{pmatrix}
1&0&0\\
0&0&0\\
0&0&-1
\end{pmatrix}
\]
 which corresponds to $a=1$ and $A=0$, that is, $h\in Z(\fr{g}_{0})$. 
Note  also that, $\fr{g}_{-}:=\fr{g}_{-2}\oplus\fr{g}_{-1}$ is a Heisenberg algebra, while $(e, h, f)$ is an $\fr{sl}_{2}$-triple, both naturally associated  to the contact grading. On the other hand, the subalgebra $\fr{p}:=\fr{g}_{\geq}=\fr{g}_{0}\oplus\fr{g}_{1}\oplus\fr{g}_{2}$ is a parabolic subalgebra of $\fr{so}^*(2n+4)$ and we will denote by 
 $P\subset\SO^*(2n+4)$  the corresponding parabolic subgroup of $\SO^*(2n+4)$.    It is well-known  that the $(4n+1)$-dimensional   real flag manifold  
  $N:=\SO^\ast(2n+4)/P$
admits an invariant contact distribution $\mathscr{D}$  (see \cite[Section 4.2]{CS09}), 
\[
\mathscr{D}\cong \SO^*(2n+4)\times_{P}\big((\fr{g}_{-1}\oplus\fr{p})/\fr{p}\big)\,.
\]
Then, for a sufficiently small open subset  $U\subset N$ of $N$ and an 1-dimensional subalgebra $\frak{t}\subset \frak{so}^*(2n+4)$  such that the corresponding right-invariant vector fields are transversal to $\mathscr{D}$ everywhere on $U$, there is a local leaf space 
\[
M:=\exp(\frak{t})\backslash U
\]
 such that $q: U\to M$ is a line bundle of scalar 2-forms on $M$. Indeed, the quaternionic structure  on $M$   induced by the isomorphism  $(\dd q)_{x}(\mathscr{D}_x)\cong T_{q(x)}M$, 
   is independent of the choice of $x$ in the fiber of $q$. Moreover,    local sections $M\to U$ correspond (via the Levi-bracket) to local almost qs-H structures on $M$ of torsion type $\mc{X}_4$. In particular, there are sections $M\to U$ corresponding to local 
   (torsion-free) qs-H structures on $M$ described in \cite{CahS}, which form one-parameter families of (local) qs-H structures, where the corresponding scalar 2-forms differ by a positive scalar, see also \cite[Remark 1.6.]{CGWPartII}. Therefore, two subalgebras $\frak{t}$ and $\hat{\frak{t}}$ provide one-parameter families of locally qs-H isomorphic  qs-H structures on $M$ if and only if they are conjugated by an element of $\SO^\ast(2n+4)$.

We will now show that  this construction can provide  
the full   classification of  all qs-H tiLas.
First we will establish an important bijection. 

\bp\label{general_approach1}
There is a bijection between one-parameter families of isomorphism classes of qs-H tiLas $(\frak{g},\sigma,Q_o, e^{2t}\omega_o)$ $(t\in\R)$ and  $\SO^\ast(2n+4)$-conjugacy classes of 1-dimensional subalgebras $\frak{t}\subset \frak{so}^*(2n+4)$ represented by elements of the form
\begin{align}\label{canelement}
\tau=\begin{pmatrix}
a & 0& d\\
0 & A & 0\\
1&0 & a
\end{pmatrix},
\end{align}
such that $a\in \sp(1)$, $A\in \fr{so}^*(2n)$ and $d\in \R$ satisfy the matrix equation
\begin{align}\label{symtest}
Xd+2AXa-Xa^2-A^2X=0
\end{align}
for all $X\in \Hn^n$.  For such a representative $\tau$, the reductive complement $\fr{m}$ in the canonical decomposition of $\fr{g}$ is given by
\[
\fr{m}=\left\{\begin{pmatrix}
0 & (AX-Xa)^*\mathbb{j}& 0\\
X & 0 & AX-Xa\\
0&-X^*\mathbb{j} & 0 
\end{pmatrix}  \right\}.
\]
\ep
\pr
We start by proving that the one-parameter family of qs-H structures $(Q, e^{2t}\om) $  $(t\in\R)$ on the local leaf space $M=\exp(\frak{t})\backslash U$ of $N$ provided by $\frak{t}\subset \frak{so}^*(2n+4)$,  with the claimed properties, gives rise to a qs-H symmetric space $(M, Q,  e^{2t}\om, s)$.   Recall that the Lie algebra of the qs-H vector fields of  $(Q, e^{2t}\om) $ is isomorphic to the quotient $N_{\frak{so}^*(2n+4)}(\frak{t})/\frak{t}$, see \cite{CapII}.  Consider now elements $\tilde X$ in $\frak{so}^*(2n+4)$ of the form
\begin{align*}
\begin{pmatrix}
0 & (AX-Xa)^*\mathbb{j}& 0\\
X & 0 & AX-Xa\\
0&-X^*\mathbb{j} & 0 
\end{pmatrix}\,,
\end{align*}
where $X\in \Hn^n$ and $A,a$ are given by the claimed representative $\tau$ of $\frak{t}$, see (\ref{canelement}). 
Then,   under our assumptions we get that
\[
[\tilde X, \tau]=\begin{pmatrix}
0 & (Xd+2AXa-Xa^2-A^2X)^*\mathbb{j}  & 0\\
0 & 0 & Xd+2AXa-Xa^2-A^2X\\
0& 0 & 0 
\end{pmatrix}=0\,.
\]
 Therefore $\tilde X\in Z_{\frak{so}^*(2n+4)}(\frak{t})/\frak{t}\subset N_{\frak{so}^*(2n+4)}(\frak{t})/\frak{t}$ for all $X\in \Hn^n$. Denote by 
 \[
 \frak{m}\subset Z_{\frak{so}^*(2n+4)}(\frak{t})/\frak{t}
 \]
 the space spanned by all such $\tilde X$. Then $\frak{g}=\frak{m}\oplus [\frak{m},\frak{m}]$ is a tiLa with involution defined as $-\id_{\frak{m}}$ in the $\fr{m}$-part and by $\id_{ [\frak{m},\frak{m}]}$ in the $[\fr{m}, \fr{m}]$-part. This is because of the inclusion $[[\frak{m},\frak{m}],\frak{m}]\subset \frak{m}$, which follows  as there is no element of $\frak{so}^*(2n+4)_{1}$ in $Z_{\frak{so}^*(2n+4)}(\frak{t})/\frak{t}.$ Then, by  construction, this is a qs-H tiLa. Indeed, we have $\omega_o\in \Lambda^2 \frak{m}^*\otimes \frak{t}$, which  corresponds  to the Levi-bracket and  is preserved by the symmetries (because $[\frak{m},\frak{m}]\subset Z_{\frak{so}^*(2n+4)}(\frak{t})/\frak{t}$).

\noindent Let us now consider two abelian subalgebras  $\frak{t},\hat{\frak{t}}\subset \frak{so}^*(2n+4)$, represented by $a,A,d$ and $\hat{a},\hat{A},\hat{d}$, respectively. Suppose that they are   conjugated by an element $g\in \SO^\ast(2n+4)$. In fact, this element has to belong in $\Sl(2,\R)\times \Sp(1)\SO^\ast(2n)$, where $\Sl(2,\R)$ corresponds to the $\fr{sl}_{2}$-triple $(e, h, f)$.  This is because  otherwise the representatives given by $a,A,d$ and $\hat{a},\hat{A},\hat{d}$ would not be mapped onto each other. Note that  the action of $\Sl(2,\R)$ maps $a,A,d$ (after rescaling)  to $e^{2t}a, e^{2t}A, e^{4t}d$, for some $t\in \mathbb{R}$.  This induces a Lie algebra isomorphism $t\id_{\frak{m}}\oplus \id_{\frak{l}}$ between  the corresponding tiLas which rescales the scalar 2-form, i.e.,  preserves the one-parameter family.    On the other hand, the action of $(b,B)\in \Sp(1)\SO^\ast(2n)$ maps $a,A,d$ to $bab^{-1},BAB^{-1},d$, and this induces    an isomorphism of qs-H tiLas  via the adjoint action of $\Sp(1)\SO^\ast(2n)$.  Thus  we conclude that
\begin{equation}\label{conjclass}
\hat{a}=e^{2t}bab^{-1}\,,\quad \hat{A}=e^{2t}BAB^{-1}\,,\quad \hat{d}=e^{4t}d
\end{equation}
 and the corresponding qs-H tiLas  belong to the same one-parameter family.  
 
\noindent Conversely,  for any almost qs-H manifold $(M, Q,\omega)$ of dimension $4n$  there is a construction given in \cite[Theorem 4]{CapII} that provides a parabolic contact structure of type $(\SO^\ast(2n+4), P)$, which  contains a line bundle $q: U\to M$ over $M$ and a contact distribution $\mathcal{D}$ on total space $U$.   This construction depends only on the conformal class of $\omega$,  hence the result will be the same for the whole one-parameter family $e^{2t}\omega$ of scalar 2-forms.   Moreover, by \cite[Theorem 5]{CapII} it follows that every qs-H transformation lifts to an isomorphism of the parabolic contact structures. In particular, for a qs-H symmetric space $(M,Q, e^{2t}\omega, s)$ in the one-parameter family and every $y\in U$, there exists a unique automorphism $s_y$ of the parabolic contact structure such that $s_y(y)=y$ and $s_y$ is a lift of $s_{q(y)}$, i.e., 
\begin{equation}\label{harmonic_d}
(\dd s_y)_{y}|_{\mc{D}_{y}}=-\id_{\mathcal{D}_y}\,.
\end{equation}
At this point, recall  that the fundamental invariant of the parabolic contact structure is  the \textsf{harmonic torsion} (in terms of \cite[p.~419]{CS09}) which is a special section of the vector bundle $\Lambda^2 \mathcal{D}^*\otimes \mathcal{D}$. Since  by (\ref{harmonic_d}) the differential $(\dd s_y)_{y}|_{\mc{D}_{y}}$  acts as $-\id$ on  $\Lambda^2  \mathcal{D}^*\otimes \mathcal{D}$, this harmonic torsion should vanish.   Therefore, the parabolic contact structure is locally isomorphic to the homogeneous space $N=\SO^\ast(2n+4)/P$, that is, $\mc{D}\cong\mathscr{D}$.  Then, by construction \cite[Theorem 4]{CapII},   there is an abelian subalgebra $\frak{t}\subset \frak{so}^*(2n+4)$ such that the corresponding right-invariant vector fields are vertical, i.e.,  belong to $\Ker(\dd q)$.  Moreover, the qs-H tiLa $(\frak{g},\sigma,Q_o, e^{2t}\omega_o)$ corresponding to the qs-H symmetric space $(M,Q, e^{2t}\omega,s)$ lifts to a subalgebra $\frak{g}'\subset N_{\frak{so}^*(2n+4)}(\frak{t})$.  Then $\sigma$ corresponds to the involution induced by   conjugation by  the elemeny $-\id\in \SO^\ast(2n)\subset P.$ Therefore, after conjugation  by an element of the subgroup $\exp(\frak{so}^*(2n+4)_2)\subset P$,  we can bring $\frak{t}$ into a form represented by $a,A,d$ as in (\ref{canelement}). Moreover,  there is a lift  $\frak{m}' \subset N_{\frak{so}^*(2n+4)}(\frak{t})$ of the reductive complement $\frak{m}$ that takes the  form
\begin{align*}
\begin{pmatrix}
0 & f(X)^*\mathbb{j}& 0\\
X & 0 & f(X)\\
0&-X^*\mathbb{j} & 0
\end{pmatrix}\,,
\end{align*}
for  all $X\in \Hn^n$ and some linear map $f : \Hn^n\to \Hn^n$.  This is because $\fr{so}^*(2n+4)_{\pm 1}$ is the $(-1)$-eigenspace of the corresponding involution. Then, a direct computation shows that the condition   $[\frak{m}',\frak{t}]=0$  implies that $f(X)=AX-Xa$ and that $Xd+2AXa-Xa^2-A^2X=0$, for all $X\in \Hn^n.$

\noindent  Recall now that by Theorem \ref{iLas_thm}  an isomorphism of  qs-H tiLas corresponds to a qs-H isomorphism between the corresponding symmetric spaces. Moreover,  by \cite[Theorem 5]{CapII} and \cite[Proposition 1.5.2]{CS09} this isomorphism should lift to left multiplication by an element of $\SO^\ast(2n+4)$.  This shows that  the corresponding  abelian subalgebra $\frak{t}$ should belong to the same $\SO^\ast(2n+4)$-conjugacy class and completes the proof. 
\pro

Let us now use Proposition \ref{general_approach1} to construct examples of qs-H tiLas.
As we will see later, these examples exhaust all possibilities and yield the full classification 
of qs-H tiLas. Next we will denote by $I_{p,q}$ the diagonal matrix with $p$ ones and $q$ negative ones on the diagonal.

\bex\label{clas_exa_1}
Let  $n$ be even and set $\mathbb{j}=\begin{pmatrix}
 0 &  I_{\frac{n}2,0} \\
 I_{0,\frac{n}2} & 0\\
\end{pmatrix}$, that is, we identify $\Hn^n$   with the linear model $[\E\Hh]$ using  a quaternionic Darboux basis, see \cite[Section A.2]{CGWPartI}.
 Let also $p, q$ be integer parameters such that  $0\leq p,q,p+q \leq \frac{n}2$ and consider the generators $\tau$  of $\fr{t}$ of the form
\begin{align*}
\tau=\begin{pmatrix}
0& 0& 0&0& 0&0&0&0\\
0& 0& 0&0& 0&0&0&0\\
0& 0& 0&0& 0&I_{q,0}&0&0\\
0& 0& 0&0& 0&0&0&0\\
0&  I_{p,0}& 0&0& 0&0&0&0\\
0& 0& 0&0& 0&0&0&0\\
0& 0& 0&0& 0&0&0&0\\
1& 0& 0&0& 0&0&0&0
\end{pmatrix}\,,
\end{align*}
 that is, $\tau$ admits a block decomposition with blocks of sizes $1,p,q,\frac{n}2-(p+q),p,q,\frac{n}2-(p+q),1$ respectively. Hence $\tau$ is a $(n+2)\times (n+2)$ quaternionic matrix.
 Then $\tau$ satisfies (\ref{symtest})  and hence by Proposition \ref{general_approach1} we obtain a qs-H tiLa   where the Lie algebra $\fr{g}$ is defined by $\fr{g}=\big([\fr{m}, \fr{m}]/\fr{t}\big)\oplus\fr{m}$, where $\fr{m}$ has the block matrix presentation
 \[
 \fr{m}=\begin{pmatrix}
0 & -X_1^*& 0&0&  0&X_2^*&0&0\\

X_1& 0 & 0&0&  0& 0 &0&0\\

Y_1&  0 & 0 & 0& 0 & 0 & 0 &X_2\\

Z_1& 0 & 0&0&  0& 0&0&0\\

Y_2& 0 & 0 & 0 & 0 & 0 & 0 &X_1\\

X_2&0 & 0&0 &0& 0&0&0\\

Z_2& 0& 0&0& 0&0&0&0\\

0 & Y_2^*& X_2^*&Z_2^* &-X_1^*&-Y_1^*&-Z_1^*&0
\end{pmatrix}.
 \]
Note that the blocks of $\fr{m}$ have  the same size with the blocks of the representative $\tau$. 
It follows that $\fr{g}$   admits the  block matrix representation
 \[
\fr{g}=\begin{pmatrix}
S_1& -X_1^*& 0&0&  0&X_2^*&0&0\\

X_1& S_2& 0&0&  0&S_3^*&0&0\\

Y_1& H_1& S_4& R_1 & S_3 & H_2 & R_2^*&X_2\\

Z_1& R_3& 0&0&  0& R_{2}&0&0\\

Y_2& H_3& S_3^*& R_4^* & S_2& -H_1^*&-R_3^*&X_1\\

X_2& S_3& 0&0 &0& S_4&0&0\\

Z_2& R_4& 0&0& 0&-R_1^*&0&0\\

H_4 & Y_2^*& X_2^*&Z_2^* &-X_1^*&-Y_1^*&-Z_1^*&S_1
\end{pmatrix}/\fr{t} \,.
\]
Here, we have $X_a, Y_a, Z_a\in\fr{m}$ with $1\leq a\leq 2$ and $S_b, H_b, R_b\in[\fr{m}, \fr{m}]$, respectively, for $1\leq b\leq 4$, with 
\[
S_1=-S_1^*,\quad S_2=-S_2^*, \quad S_4=-S_4^*,\quad 
H_2=H_2^*,\quad H_3=H_3^*, \quad H_4=H_4^*\,.
\]
Then we see that the matrices  $X_1, X_2, S_1, \ldots, S_4$ span the semisimple part of $\fr{g}$,
which is isomorphic to $\sp(p+1, q)$. Moreover, the radical $\fr{r}\subset\fr{g}$ of $\fr{g}$ decomposes into two distinguished parts
\[
\fr{r}=\fr{r}_{\rm qH}/\fr{t}\oplus\fr{r}_{\rm deg}\,,
\]
where $\fr{r}_{\rm qH}$ is generated by   the matrices $Y_1, Y_2, H_1, \ldots, H_4$  and  $\fr{r}_{\rm deg}$ is generated by 
$Z_1, Z_2, R_1, \ldots, R_4$. 
Note that   $\fr{r}_{\rm qH}$ can be arranged into  a block matrix decomposition, of sizes $1, p, q$, i.e., 
\[
\begin{pmatrix}
H_4 & Y_2^* & -Y_1^*\\
Y_1& H_1 & H_2\\
Y_2& H_3 & -H_1^*
\end{pmatrix}
\]
and, by dimensional reasons,  is isomorphic to the space of quaternionic Hermitian scalar products on the standard representation of $\fr{sp}(p+1, q)$. In particular, $\tau$ belongs to $\fr{r}_{\rm qH}$  and represents  the non-degenerate quaternionic Hermitian scalar product defining the Lie algebra $\fr{sp}(p+1, q)$.  We mention that:\\
\noindent $\bullet$ For  $p+q=0$, we have $\tau=H_4$ and the Lie algebra $\fr{g}$ is exactly the linear model $[\E\Hh]$.\\
\noindent $\bullet$ For general $p, q$   the orbit of the semisimple part $\fr{sp}(p+1, q)$ is  the quaternionic pseudo-K\"ahler symmetric space $\Sp(p+1, q)/(\Sp(1)\Sp(p, q))$. \\
\noindent Let us now highlight a few low-dimensional cases, where special isomorphisms of simple Lie algebras occur. 
\begin{itemize}
\item[(1)]  Let us consider the case $n=2$, $p=1$ and $q=0$.  In this case, the  block sizes are $(1,1,0,0,1,0,0,1)$ and Lie algebra $\fr{g}$  has the matrix presentation
\[
\fr{g}=\left\{
\begin{pmatrix}
S_1 & -X_1^* & 0 & 0\\
X_1 & S_2 & 0 & 0\\
Y_2 & H_3 & S_2 & X_1\\
H_4 & Y_2^* & -X_1^* & S_1
\end{pmatrix} : S_1, S_2\in\fr{sp}(1), X_1, Y_2\in\Hn, H_3, H_4\in\R
\right\}/\fr{t}\,,
\]
thus $\fr{g}$ is 15-dimensional. 
Using the isomorphism $\fr{sp}(2)\cong\fr{so}(5)$ we may express the Levi decomposition as  $\fr{g}=\fr{so}(5)\oplus\R^5$, where $\fr{so}(5)$ acts on $\R^{5}$ by the standard representation.  Moreover, we see that  $\fr{l}=\fr{so}(4)\oplus\R$, and the reductive complement $\fr{m}$ is 8-dimensional,
\[
\fr{m}=\left\{
\begin{pmatrix}
0 & -X_1^* & 0 & 0\\
X_1 & 0 & 0 & 0\\
Y_2 & 0 & 0 & X_1\\
0 & Y_2^* & -X_1^* & 0
\end{pmatrix} :   X_1, Y_2\in\Hn\right\}\,.
\]
This yields the  non-semisimple symmetric space  
\[
(\SO(5)\ltimes\R^5)/(\SO(4)\times\R)\,.
\]
 In this case,  the semisimple quaternionic K\"ahler symmetric space is the coset $\SO(5)/\SO(4)$, that is, the 4-sphere $\Ss^4\cong\Hn{\sf P}^1$. 
\item[(2)] Let us now consider the case $n=2$, $p=0$ and $q=1$. Here the  sizes of the blocks are $(1,0,1,0,0,1,0,1)$ and Lie algebra $\fr{g}$  has the matrix presentation
\[
\fr{g}=\left\{\begin{pmatrix}
S_1& 0 & X_2^*& 0\\
Y_1& S_4&  H_2 &X_2\\
X_2&  0&  S_4& 0 \\
H_4 &  X_2^*& -Y_1^* &S_1
\end{pmatrix} : S_1, S_4\in\fr{sp}(1), X_2, Y_1\in\Hn, H_2, H_4\in\R
\right\}/\fr{t} \,.
\]
Hence $\fr{g}$ is again 15-dimensional. Using   the isomorphism $\fr{sp}(1, 1)\cong\fr{so}(4, 1)$, we can express the Levi decomposition of $\fr{g}$ as $\fr{g}=\fr{so}(4, 1)\oplus\R^5$, where $\R^5$ is the standard representation of $\fr{so}(4, 1)$. 
In this case we get $\fr{l}=\fr{so}(4)\oplus\R$ and 
\[
\fr{m}=\left\{\begin{pmatrix}
0& 0 & X_2^*& 0\\
Y_1& 0&  0 &X_2\\
X_2&  0&  0& 0 \\
0 &  X_2^*& -Y_1^* & 0
\end{pmatrix} :   X_2, Y_1\in\Hn\right\} \,.
\]
This corresponds to the  8-dimensional non-semisimple symmetric space  
\[
(\SO(4, 1)\ltimes\R^5)/(\SO(4)\times\R)
\]
 where the semisimple quaternionic K\"ahler symmetric space    is the non-compact dual of the previous one, that is,  the symmetric space $\SO(4, 1)/\SO(4)$.
\end{itemize}
\eex

\bex\label{clas_exa_2}
Suppose that $n$ is odd and set 
\[
\mathbb{j}=\begin{pmatrix}
 0 & 0& I_{\frac{n-1}2,0} \\
 0& j& 0\\
 I_{0,\frac{n-1}2} & 0& 0
\end{pmatrix}\,.
\]
 This choice identifies  $\Hn^n$ with the linear model $[\E\Hh]$ by a generalization of  a quaternionic Darboux basis, see \cite[Proposition A.9]{CGWPartI}.
Let $p, q$ be integer parameters such that $0\leq p,q,p+q \leq \frac{n-1}2$,
and consider the generator $\tau\in\fr{t}$ of the form
\begin{align*}
\tau=\begin{pmatrix}
0& 0& 0&0& 0&0&0&0&0\\
0& 0& 0&0& 0&0&0&0&0\\
0& 0& 0&0& 0&0&I_{q,0}&0&0\\
0& 0& 0&0& 0&0&0&0&0\\
0& 0& 0&0& 0&0&0&0&0\\
0&  I_{p,0}& 0&0& 0&0&0&0&0\\
0& 0& 0&0& 0&0&0&0&0\\
0& 0& 0&0& 0&0&0&0&0\\
1& 0& 0&0& 0&0&0&0&0
\end{pmatrix}.
\end{align*}
Hence,  in this case the $(n+2)\times(n+2)$ quaternionic matrix $\tau$ admits a block decomposition   with blocks of sizes  $1,p,q,\frac{n-1}2-(p+q),1,p,q,\frac{n-1}2-(p+q),1$.  It is easy to see that 
 $\tau$ satisfies (\ref{symtest}),  hence by Proposition \ref{general_approach1} we obtain a qs-H tiLa  with Lie algebra $\fr{g}$ defined by $\fr{g}=\big([\fr{m}, \fr{m}]/\fr{t}\big)\oplus\fr{m}$,
 where the reductive complement $\fr{m}$ is given by
  \[
\fr{m}=\begin{pmatrix}
0& -X_1^*& 0&0&  0&0&X_2^*&0&0\\

X_1&0& 0&0&  0&0& 0&0&0\\

Y_1& 0& 0& 0 & 0& 0 & 0 & 0&X_2\\

Z_1& 0& 0&0&  0& 0& 0&0&0\\

Z_3& 0& 0& 0& 0&0 & 0& 0& 0\\

Y_2& 0& 0& 0 &0 & 0&0&0&X_1\\

X_2&0& 0&0 &0 & 0& 0&0&0\\

Z_2& 0& 0&0& 0& 0& 0&0&0\\

0 & Y_2^*& X_2^*& Z_2^* & -Z_3^*j & -X_1^*&-Y_1^*&-Z_1^*&0
\end{pmatrix}.
\]
Note that  the  block structure of $\fr{m}$ matches the block sizes   of the representative $\tau$. 
It follows that $\fr{g}$ admits the following block matrix representation
 \[
\fr{g}=\begin{pmatrix}
S_1& -X_1^*& 0&0&  0&0&X_2^*&0&0\\

X_1& S_2& 0&0&  0&0& S_3^*&0&0\\

Y_1& H_1& S_4& R_1 & R_6& S_3 & H_2 & R_2^*&X_2\\

Z_1& R_3& 0&0&  0& 0& R_{2}&0&0\\

Z_3& R_5& 0& 0& 0&0 & jR_6^*& 0& 0\\

Y_2& H_3& S_3^*& R_4^* &-R_5^*j & S_2& -H_1^*&-R_3^*&X_1\\

X_2& S_3& 0&0 &0 & 0& S_4&0&0\\

Z_2& R_4& 0&0& 0& 0& -R_1^*&0&0\\

H_4 & Y_2^*& X_2^*& Z_2^* & -Z_3^*j & -X_1^*&-Y_1^*&-Z_1^*&S_1
\end{pmatrix}/\fr{t} 
\]
with $X_1, X_2, Y_1, Y_2,  Z_1, \ldots, Z_3\in\fr{m}$  and $S_b, H_b, R_1,\ldots, R_6\in[\fr{m}, \fr{m}]$, respectively, for $1\leq b\leq 4$, with 
\[
S_1=-S_1^*,\quad S_2=-S_2^*, \quad S_4=-S_4^*,\quad 
H_2=H_2^*,\quad H_3=H_3^*, \quad H_4=H_4^*\,.
\]
Similarly with the previous example we see that  that the matrices  $X_1, X_2, S_1, \ldots, S_4$ span the semisimple part of $\fr{g}$,
which is isomorphic to $\sp(p+1, q)$. Moreover,  in this case the decomposition $\fr{r}=\fr{r}_{\rm qH}/\fr{t}\oplus\fr{r}_{\rm deg}$ of the   radical $\fr{r}\subset\fr{g}$ of $\fr{g}$ is such that
$\fr{r}_{\rm qH}={\rm span}\{Y_1, Y_2, H_1, \ldots, H_4\}$ and 
$\fr{r}_{\rm deg}={\rm span}\{Z_1, Z_2, Z_3, R_1, \ldots, R_6\}$, respectively. 
Again the component $\fr{r}_{\rm qH}$ is isomorphic to the space of quaternionic Hermitian scalar products on the standard representation of $\fr{sp}(p+1, q)$, where $\tau\in\fr{r}_{\rm qH}$ represents  the non-degenerate quaternionic Hermitian scalar product defining the Lie algebra $\fr{sp}(p+1, q)$. In particular, 
for  $p+q=0$, we similarly obtain $\tau=H_4$ and the Lie algebra $\fr{g}$ is exactly the linear model $[\E\Hh]$,
while   for general $p, q$  the orbit of the semisimple part $\fr{sp}(p+1, q)$ is once more the quaternionic pseudo-K\"ahler symmetric space $\Sp(p+1, q)/(\Sp(1)\Sp(p, q))$. \\
\noindent Let us discuss a  low-dimensional case.  Suppose that $n=3$, $p=1$, $q=0$. In this case $\fr{g}$ has the matrix presentation  
\[
\fr{g}=\left\{\begin{pmatrix}
S_1& -X_1^*  &  0&0 &0\\

X_1& S_2 &  0&0 &0\\

Z_3& R_5 & 0&0  & 0\\

Y_2& H_3 &-R_5^*j & S_2 &X_1\\

-H_3 & Y_2^* & -Z_3^*j & -X_1^*& S_1
\end{pmatrix} : X_1, Y_2, Z_3, R_5\in\Hn, S_1, S_2\in\fr{sp}(1), H_3\in\R \right\}
\]
and we obtain a 12-dimensional non-semisimple qs-H symmetric space with
  reductive complement $\fr{m}$ having the form
\[
\fr{m}=\left\{\begin{pmatrix}
0& -X_1^*  &  0&0 &0\\

X_1&0 &  0&0 &0\\

Z_3& 0 & 0&0  & 0\\

Y_2& 0 & 0 & 0 &X_1\\

0 & Y_2^* & -Z_3^*j & -X_1^*&0 
\end{pmatrix} : X_1, Y_2, Z_3\in\Hn\right\}\,.
\]
Note that  the block $\begin{pmatrix} S_1 & -X_1^*\\ X_1	 & S_2\end{pmatrix}$ represents  the semisimple part, which corresponds to the  quaternionic K\"ahler symmetric space $\Hn{\sf P}^{1}=\Sp(2)/(\Sp(1)\times\Sp(1))$. The solvable part of $\fr{g}$,  as a vector space, is the direct sum of the standard representation $\Hn^2$ of $\Sp(2)$ and the space of the quaternionic Hermitian scalar products on $\Hn^2$. However,  observe that the Lie bracket on the solvable part contains the pairing 
\[
[(u, v), (w, z)]=\begin{pmatrix}
z^*ju-v^*jw & z^*jv-w^*jz\\
 w^*ju-u^*jw & w^*jv-u^*jz
\end{pmatrix}/(\mu\Id)\,,\quad\text{for}\quad (u, v), (w, z)\in\Hn^2, \mu\in\R\,,
\]
where $\mu\Id$ represents $\fr{t}$, hence  in this case the radical is not abelian.
 \eex

 Let us now describe the realization of the semisimple qs-H tiLas (see Theorem \ref{ClasThm}) via the approach that encodes  Proposition \ref{general_approach1}.

\bex\label{clas_exa_3}
For general $n>1$ consider the $(n+2)\times(n+2)$ quaternionic matrix 
\[
\tau=\begin{pmatrix}
i &0& -1\\
0 & 0 &0\\
1&0 & i
\end{pmatrix}
\]
with a block decomposition  into blocks of sizes $1, n, 1$.
Setting  $\mathbb{j}=j$ we can identify $\Hn^n$  with the linear model $[\E\Hh]$ by means of a skew-Hermitian basis, see \cite[Section A.2]{CGWPartI}. 
It is easy to see that $\tau$ satisfies (\ref{symtest}), hence by Proposition \ref{general_approach1} we obtain a qs-H tiLa  with Lie algebra $\fr{g}$ defined by $\fr{g}=\big([\fr{m}, \fr{m}]/\fr{t}\big)\oplus\fr{m}$,
 where this time the reductive complement $\fr{m}$ is given by 
  \[
\fr{m}=\left\{\begin{pmatrix}
0 &   iX^*j & 0         \\
X & 0 &  -Xi   \\
0 & -X^*j & 0 
\end{pmatrix}  : X\in \Hn^n\right\}.
\]
   Again the block structure of $\fr{m}$ matches   the   block sizes   of $\tau$. 
It follows that $\fr{g}$ admits the following block matrix representation
 \[
\fr{g}=\left\{\begin{pmatrix}
S_1 &   iX^*j & 0         \\
X & S_2 &  -Xi   \\
0 & -X^*j & S_1 
\end{pmatrix} : S_1\in\fr{u}(1), S_2\in\fr{so}^*(2n), X\in\Hn^n\right\},
\]
which, obviously,   corresponds to the simple qs-H symmetric space $\SO^*(2n+2)/\SO^*(2n)\U(1).$
\eex

\bex\label{clas_exa_4}
For $n=p+q>1$ consider the $(n+2)\times(n+2)$ quaternionic matrix 
\[
\tau=\begin{pmatrix}
0 &0& -1\\
0 & jI_{p,q} &0\\
1&0 & 0
\end{pmatrix}
\]
with a block decomposition  into blocks of sizes $1, n, 1$.
Moreover, set $\mathbb{j}=j$ such that we can identify   $\Hn^n$ with $[\E\Hh]$ by means of a skew-Hermitian basis. 
The matrix $\tau$ satisfies (\ref{symtest}), hence once more by Proposition \ref{general_approach1} we obtain a qs-H tiLa  with Lie algebra $\fr{g}$ defined by $\fr{g}=\big([\fr{m}, \fr{m}]/\fr{t}\big)\oplus\fr{m}$. In this case 
the  reductive complement $\fr{m}$ is given by 
  \[
\fr{m}=\left\{\begin{pmatrix}
0 &    X^*I_{p, q} & 0         \\
X & 0 &  jI_{p, q}X   \\
0 & -X^*j & 0 
\end{pmatrix} : X\in \Hn^n\right\},
\]
 where the blocks   have  the same size with the blocks  of $\tau$. 
It follows that $\fr{g}$ admits the following block matrix representation
 \[
\fr{g}=\left\{\begin{pmatrix}
S_1 &  X^*I_{p, q} & 0         \\
X & S_2 &  jI_{p, q}X  \\
0 & -X^*j & S_1 
\end{pmatrix} : S_1\in\fr{sp}(1), S_2\in\fr{u}(p, q), X\in\Hn^n\right\}\,.
\]
This corresponds to the simple qs-H symmetric space $\SU(2n+2)/(\SU(2)\SU(p,q)\U(1)).$
\eex

\bex\label{clas_exa_5}
For any even $n$ consider the $(n+2)\times(n+2)$ quaternionic matrix $\tau$
\[
\tau=\begin{pmatrix}
0 &0&0& 1\\
0 & I_{\frac{n}2,0} &0& 0\\
0 &0& I_{0,\frac{n}2} & 0\\
1&0 & 0 & 0
\end{pmatrix}
\]
 with block decomposition of sizes $1,\frac{n}2,\frac{n}2,1$.
Moreover,   identify $\Hn^n$ with $[\E\Hh]$ by means of a quaternionic Darboux basis, by setting 
\[
\mathbb{j}=\begin{pmatrix}
 0 &  I_{\frac{n}2,0} \\
 I_{0,\frac{n}2} & 0\\
\end{pmatrix}.
\]
Again, the matrix $\tau$ satisfies (\ref{symtest}) and  Proposition \ref{general_approach1} gives  a qs-H tiLa  with Lie algebra $\fr{g}$ defined by $\fr{g}=\big([\fr{m}, \fr{m}]/\fr{t}\big)\oplus\fr{m}$. In this case 
the  reductive complement $\fr{m}$ is given by 
  \[
\fr{m}=\left\{\begin{pmatrix}
0 &    X_2^* & X_1^*  & 0       \\
X_1 & 0 &  0 & X_1   \\
X_2 & 0 & 0 & -X_2 \\
0 & X_2^* & -X_1^* & 0 
\end{pmatrix} : X_1, X_2\in\Hn^{\frac{n}{2}}\right\},
\]
where  the blocks   have  the same size with the blocks  of $\tau$. 
Therefore, for the Lie algebra $\fr{g}$  we obtain the  block matrix representation
 \[
\fr{g}=\left\{\begin{pmatrix}
S_1 &    X_2^* & X_1^*  & 0       \\
X_1 & S_2 &  0 & X_1   \\
X_2 & 0 & -S_2^* & -X_2 \\
0 & X_2^* & -X_1^* & S_1
\end{pmatrix} : S_1\in\fr{sp}(1), S_2\in\fr{gl}(\frac{n}{2}, \Hn), X_1, X_2\in\Hn^{\frac{n}{2}}\right\}\,.
\]
This corresponds to the simple symmetric space $\Sl(\frac{n}2+1,\mathbb{H})/(\Gl(1,\mathbb{H})\Sl(\frac{n}2,\mathbb{H}))$.
\eex

We will now prove that Examples \ref{clas_exa_1}-\ref{clas_exa_5} provide the full classification 
of  qs-H tiLas.  For this step we will use the classification of normal forms  of Jordan blocks of the Lie algebras $\fr{sp}(1)$ and $\fr{so}^*(2n)$   presented in \cite[Table III]{Djokovic}. Next we will denote by   $J_{r}(\lambda)$ the standard Jordan blocks of size $r$ corresponding to the eigenvalue $\lambda\in\Hn$, with the units located under the diagonal, i.e.,
\[
J_{r}(\lambda)=\begin{pmatrix}
 \lambda & &  & \\
 1 & \lambda & & \\
    & \ddots & \ddots & \\
    &            & 1         & \lambda
    \end{pmatrix}\,.
\]

\bt\label{main_thm_tiLa}
 There is a bijection between the  one-parameter families  of isomorphism classes of qs-H tiLas $(\fr{g},\sigma, Q_o, e^{2t}\om_o)$ $(t\in \R)$, and the quaternionic matrices $\tau$ presented
 in   Examples {\rm \ref{clas_exa_1}-\ref{clas_exa_5}}.
 \et
\pr
For an abelian subalgebra $\fr{t}\subset\fr{so}^*(4n+2)$  we can pick a representative  $\tau\in\fr{t}$ of the form \eqref{canelement} under the action defined in \eqref{conjclass}, and test whether   it satisfies the condition given in \eqref{symtest}. Recall that $\tau$ depends on $a\in\fr{sp}(1)$, $A\in\fr{so}^*(2n)$ and $d\in\R$. 

According to \cite{Djokovic} and the classification of  Jordan blocks for $\fr{so}^*(2n)$,  and also for $\fr{sp}(1)$,  we can construct the matrix $\tau$ by combining the following normal forms:  
 \begin{enumerate} 
\item $A\in\fr{so}^*(2n)$ is in block diagonal form with blocks (see cases (29)-(31) in \cite[Table III]{Djokovic})
\[
J_{n_u}(0)\,,\qquad J_{n_v}(b_vj)\,, \qquad J_{n_w}(\beta_w)\oplus J_{n_w}(-\beta_w)^*\,,
\] 
for $u\in U,v\in V,w\in W$ indexing the blocks, with 
\[
2n=\sum_{u\in U}n_u+\sum_{v\in V}n_v+2\sum_{w\in W}n_w\,.
\]
Moreover, we have $0<b_v\in \R$ and $\beta_w\in \C$ with 
$\Re(\beta_w)>0,\Im(\beta_w)\geq 0$ for $v\in V, w\in W$.
 Note that some of the cases have two possibilities for $\mathbb{j}$, see cases (30)-(31)  in \cite[Table III]{Djokovic}.\footnote{We   mention  that in  \cite{Djokovic}  the quaternions are represented by $2\times 2$ complex matrices.}
 \item $a=\mu\cdot i\in\fr{sp}(1)$, for some scalar $\mu\in\R$, and  $d\in\R$. \end{enumerate} 
 \noindent Now, according to the action defined in (\ref{conjclass}) we still have the freedom to rescale  $a, A, d$ by $e^{2t} a, e^{2t}A$ and $e^{4t}d$, respectively. 
Clearly, we can investigate the condition \eqref{symtest} for each block separately; Note   that a combination of solutions of \eqref{symtest} can be glued  together to a solution if and only they have the same $a$ and $d$. 

\noindent Applying this procedure, we first observe that in \eqref{symtest}  only the  matrix $A^2$   contributes for order-three nilpotent elements. Thus we can only have  Jordan blocks of size  $n_u\leq 2, n_v\leq 2$ and  $n_w\leq 2$.  For order-two nilpotent elements, we see that if $n_u=n_v=n_w=2$ for $u\in U, v\in V,w\in W$, then the conditions 
\[
\{x a=0,\quad 2b_vj x+x a=0, \quad 2\beta_w x+x a=0\}
\]
should be satisfied for all $x\in \mathbb{H}$.  Therefore, if we have $a=0$, then $n_u=2$ is possible and $d=0$ follows from \eqref{symtest}, while a direct computation shows that the other cases are not possible, i.e., the cases $n_v=1,n_w=1$ for all $v\in V$ and $w\in W$. 

\noindent In the situation when $n_u=n_v=n_w=1$, for $u\in U, v\in V$ and $w\in W$,  by \eqref{symtest}  we obtain that the conditions 
\[
\{xd-xa^2=0, \quad xd+2b_vjxa-xa^2+(b_v)^2x=0,\quad xd+2\beta_w xa-xa^2-\beta_w^2x=0\}
\]
should be satisfied for   all $x\in \Hn$. Thus here, after rescaling, we deduce that  the pairs $(a=0,d=0)$ and $(a=i, d=-1)$ are the only solutions for $u\in U$, the triple $(a=0,b_v=1,d=-1)$ is the only solution for $v\in V$ and the triple $(a=0,d=1,\beta_w=1)$ is the only solution for $w\in W$.

\noindent All together the above analysis yields the  following possible combinations:\\
(i)  $a=d=0$ and $A$ consists of blocks $J_{1}(0)$ and $J_{2}(0)$ and  this gives rise to the representative $\tau\in\fr{t}$ described in Examples \ref{clas_exa_1} and \ref{clas_exa_2}, respectively.
   This is done by comparing the convention for $\mathbb{j}$ in \cite{Djokovic} (see the column $K$ in Table III of  \cite{Djokovic}, depending on $\kappa=\pm 1$,  for the entry (29)), with our convention of $\mathbb{j}$ in the  quaternionic-Darboux bases from \cite[Section A.2]{CGWPartI}. This means that for $\kappa=1$ the conventions are the same but for $\kappa=-1$ the two conventions differ by a conjugate transposition.  Moreover, note that  the parameters $p, q$ appearing in these two examples, can be understood in the following way: 
   \begin{itemize}
   \item The parameter $p$  is the number of the Jordan blocks $J_{2}(0)$ satisfying  $\kappa=1$.
\item The  parameter $q$ is the number of the Jordan blocks $J_{2}(0)$ satisfying  $\kappa=-1$.
\end{itemize}
 Note   that for the second case,  the Jordan block $J_{2}(0)$  appears as its conjugate transpose in the  quaternionic-Darboux basis.  \\
\noindent (ii)  $a=i,d=-1$ and $A=0$. This case yields  the representative $\tau\in\fr{t}$ described in Example  \ref{clas_exa_3}. \\
\noindent (iii)   $a=0,d=-1$ and $A$ consists of Jordan blocks of type $J_{1}(j)$. This corresponds to the representative $\tau\in\fr{t}$ described in Example  \ref{clas_exa_4}. In particular the parameters  $p,q$ are encoded exactly as in case (i), by inspecting entry (30) rather than entry (29)   in    \cite[Table III]{Djokovic}. \\
\noindent (iv)  $a=0,d=1$ and $A$ consists of blocks of type $J_{1}(1)\oplus J_{1}(-1)$.
This  gives us  the representative $\tau\in\fr{t}$ described in Example  \ref{clas_exa_5}.  

Since these cases exhaust all admissible combinations, this completes the proof.
\pro

Theorem \ref{THM2} in introduction now follows as a combination of Theorems \ref{ClasThm} and \ref{main_thm_tiLa}.

\br\label{Kill_form_m}
The restriction $B_{\fr{m}}$ of the  Killing form  of $\fr{so}^{*}(2n+4)$  to $\fr{m}$ for the qs-H tiLas from Theorem \ref{main_thm_tiLa}
is a multiple of the real trace
\begin{align*}
&\tr_{\R}\left(\begin{pmatrix}
0 & (AX-Xa)^*\mathbb{j}& 0\\
X & 0 & AX-Xa\\
0&-X^*\mathbb{j} & 0, 
\end{pmatrix}\cdot
 \begin{pmatrix}
0 & (AY-Ya)^*\mathbb{j}& 0\\
Y & 0 & AY-Ya\\
0& -Y^*\mathbb{j} & 0, 
\end{pmatrix}\right)\\
&=\Re\big((AX-Xa)^*\mathbb{j}Y-X^*\mathbb{j}(AY-Ya)\big)+\tr_{\R}\big(X(AY-Ya)^*\mathbb{j}-(AX-Xa)Y^*\mathbb{j}\big)\,.
\end{align*}
For the Examples  \ref{clas_exa_1} and \ref{clas_exa_2} we have $a=0$ and $A\in\fr{so}^*(2n)$ is not of full rank. For these cases we see that $B_{\fr{m}}$ is degenerate.  Since however $\fr{g}\subset\fr{so}^*(2n+4)$ is a subalgebra of  $\fr{so}^*(2n+4)$, the degeneracy of  $B_{\fr{m}}$ implies the degeneracy of the Killing form of $\fr{g}$. This gives an alternative explanation of the fact that the  symmetric spaces corresponding to these tiLas  are  {\it not} semisimple.
\er

 
\section{Homogeneous torsion-free qs-H structures}\label{section6}
We conclude this work by considering  the following interesting task:

\smallskip
 {\bf Question:} {\it Is it possible to adapt  the  construction described in Section \ref{sec_gen_class}  to  obtain \textsf{non-symmetric} homogeneous manifolds admitting  an invariant qs-H structure?}

\smallskip
Regarding this question, we first point out that  from
\cite[Table 4]{S01} it  is known that  the homogeneous manifolds  in our question  cannot  have full holonomy $\SO^*(2n)\Sp(1)$. 
Below we will prove that homogeneous almost qs-H manifolds which are torsion-free, i.e., homogeneous qs-H manifolds, are necessarily symmetric spaces, and thus classified by the results in Section \ref{sec_gen_class}.
 
  Therefore, let us turn  our attention again to the construction described in Section~\ref{sec_gen_class} and, this time, consider 1-dimensional abelian subalgebras 
 $\fr{t}\subset  \frak{so}^*(2n+4)$
  generated by general elements  $\hat{\tau}\in  \frak{so}^*(2n+4)$ of the form
\[ 
\hat{\tau}=\begin{pmatrix}
a & C^*\mathbb{j}& d\\
0 & A & C\\
1&0 & a
\end{pmatrix},
\]
 with $a\in \sp(1)$, $A\in \fr{so}^*(2n)$,  $d\in \R$ and $C\in \Hn^n$.
Obviously, $\hat\tau$ is  of the form $(\ref{canelement})$  if and only if $C=0$.

For $C\neq 0$ next we will examine when the qs-H structure provided by the construction
is a  non-symmetric homogeneous  qs-H manifold.
  In this case, we will show that the condition
\[
[\hat{\tau},\fr{m}]\subset \fr{t},
\]
cannot be satisfied for sufficiently large subspace $\fr{m}\subset  \frak{so}^*(2n+4)$, proving the following result, that is, Theorem \ref{THM1}.
\bt\label{homogeneous_1}
A homogeneous   qs-H manifold is a qs-H symmetric space.
\et
\pr
As we want to have elements of $\fr{m}\subset  \frak{so}^*(2n+4)$ in each tangent direction, we consider the elements
\begin{align*}
\tilde X=\begin{pmatrix}
\lambda(X)+\mu(X) & f(X)^*\mathbb{j}& g(X)\\
X & \gamma(X) & f(X)\\
0&-X^*\mathbb{j} & -\lambda(X)+\mu(X)
\end{pmatrix}\,,
\end{align*}
for  all $X\in \Hn^n$ and some linear maps $f : \Hn^n\to \Hn^n$, $g,\lambda : \Hn^n\to \R$, $\mu: \Hn^n\to \fr{sp}(1)$ and $\gamma:  \Hn^n\to \fr{so}^*(2n)$.  If $[\hat{\tau},\tilde X]\subset \fr{t}$ is satisfied for all $X\in \Hn^n$, that is $\hat{\tau}$ represents a homogeneous   qs-H structure, then the 1-form $\lambda+\mu+\gamma: \fr{m}\to \R\oplus \fr{sp}(1) \oplus \fr{so}^*(2n)$ defines the connection 1-form of a torsion-free $\SO^*(2n)\Sp(1)$-connection.  In particular, $\lambda=0$ has to be satisfied as this part is not valued in $\fr{so}^*(2n)\oplus \fr{sp}(1)$. Therefore, the grading component $-2$ of $[\hat{\tau},\tilde{X}]$ vanishes, and thus $[\hat{\tau},\tilde X]=0$. Moreover, we see that the position $(1,1)$ of the matrix $[\hat{\tau},\tilde X]$ reads as 
\[
a\mu(X)+C^*\mathbb{j}X-\mu(X)a-g(X)\,.
\]
Since $[\hat{\tau},\tilde X]=0$ we thus   get the condition 
\begin{eqnarray*}
\Im(a\mu(X)+C^*\mathbb{j}X-\mu(X)a-g(X))&=&\frac12(C^*\mathbb{j}X+X^*\mathbb{j}C)\\
&=&g_1(C,X)i+g_2(C,X)j+g_3(C,X)k=0\,,
\end{eqnarray*}
 for all $X\in \Hn^n$, where $g_{1}, g_2, g_3$ are the pseudo-Hermitian metrics on the linear model $[\E\Hh]$ (see \cite[Proposition 2.15]{CGWPartI}).   As the metrics $g_1,g_2,g_3$ are non-degenerate,  $C=0$ follows. Therefore, $\hat\tau$ is of the form (\ref{canelement}) and consequently the homogeneous qs-H manifold is a qs-H symmetric space. 
\pro

Let us emphasize that being torsion-free is a key condition in Proposition \ref{homogeneous_1}; this is because,  
in general, the output of this construction  is  a conformally qs-H structure, see \cite{CapII}. In terms of intrinsic torsion components this means that the almost  qs-H structure  obtained in this way is of type  $\mc{X}_4$ (see \cite[Remark 1.6]{CGWPartII}).   Let us indeed  present  an example of  a homogeneous almost  qs-H manifold with nontrivial torsion lying in the component $\mc{X}_4$.

\bex\label{confromal_s}
Let us consider the element 
\[
\hat{\tau}=\begin{pmatrix}
i &   -2i+2j & 0  & 3       \\
0 & 2j &  0 &2-2k   \\
0 & 0 & 0 &0 \\
1 & 0 &0 & i
\end{pmatrix}
\]
with block decomposition of blocks of size $1,1,n-1,1$ and set $\mathbb{j}=j.$ Then the subspace $\fr{m}\subset  \frak{so}^*(2n+4)$  is given by
\[
\begin{pmatrix}
X_1^*(k-1) +2(x_2-x_3)i&  2X_1^* +iX_1^*j & iX_2^*j  & 2(x_2-x_3)       \\
X_1 & 0 &  ( i-j)X_2^*j &2jX_1-X_1i   \\
X_2 &- X_2(1+k) & 0 &-X_2i \\
0 &-X_1^*j & -X_2^*j  &(1+k)X_1 +2(x_2-x_3)i
\end{pmatrix},
\]
for $X_1:=x_1+x_2i+x_3j+x_4k\in \Hn$ and $X_2\in\Hn^{n-1}$.
We can directly check
\[
[\hat{\tau},\fr{m}]=2(x_4-x_1)\hat\tau
\]
and thus  this induces a homogeneous almost qs-H structure, where the torsion belongs to the $\mc{X}_4$-component.    Recalling that  $\mc{X}_{4}\cong[\E\Hh]$ (see \cite{CGWPartI}), one  can interpret  the torsion by the  1-form $\lambda$ appearing in the proof of Theorem \ref{homogeneous_1}. This is given by  
\[
\lambda((x_1+x_2i+x_3j+x_4k,X_2)^{t})=x_4-x_1\,.
\]
  Obviously, analyzing the full qs-H automorphism algebra is a nontrivial task. We observe, however, that for $n=2$  the qs-H automorphism algebra contains a 9-dimensional solvable subalgebra generated by $\fr{m}$.
\eex

%
%


\end{document}